\documentclass[11pt,a4paper]{article}
\usepackage{mathrsfs, amsmath, amsthm, amssymb, bm, mathtools, thm-restate}
\usepackage{ulem} 
\usepackage{cancel}

\usepackage{soul}

\newtheorem{theorem}{Theorem}
\newtheorem{lem}{Lemma}
\newtheorem{cor}{Corollary}
\theoremstyle{definition}
\newtheorem{definition}{Definition}

\newtheorem{question}{Question}

\usepackage[left=25mm, top=25mm, bottom=25mm, right=25mm]{geometry}
\usepackage[T1]{fontenc}
\usepackage[inline]{enumitem}
\usepackage[square, numbers, sort&compress]{natbib}

\newcommand{\dd}{\textsf{Del}}

\newcommand{\ed}{\textsf{Del}}

\newcommand{\ds}{\textsf{Del-Save}}

\usepackage[colorinlistoftodos,prependcaption,textsize=scriptsize, textwidth=20mm,]{todonotes}
\definecolor{LightBlue}{RGB}{173,216,230}

\usepackage{tcolorbox}

\begin{document}
	\title{ Arc-weighted acyclic orientation of graphs}
		\author{Huan Zhou\thanks{School of Mathematical Sciences, Zhejiang Normal University, Jinhua, China. Emails:huanzhou@zjnu.edu.cn, xdzhu@zjnu.edu.cn.}
\and
Jialu Zhu \thanks{Data Science Institute, Shandong University, Jinan, China. Email: jialuzhu@zjnu.edu.cn.}
\and
Xuding Zhu\footnotemark[1]\thanks{Grant numbers: NSFC 12371359.}}

        \maketitle
        
	\begin{abstract}
  Let $D$ be a digraph, and let $w: E(D) \to \{1,2,\ldots \}$ be a positive integer weight assignment on the arcs of $D$.  An arc $e=(u,v)$ is called dominating if $w(e) > d_{(D,w)}^+(v)$.  We say $(D,w)$ is acyclic if every non-empty sub-digraph $(D', w)$ has a dominating arc.  
 For a graph $G$ and a mapping $f \in \mathbb{N}^G$, we say  $G$ is {arc-weighted $f$-degenerate}  if there is an arc-weighted orientation $(D, w)$ of $G$ (i.e., an orientation $D$ of $G$ together with a positive integer weight assignment $w$) such that $(D,w)$ is acyclic and  $d_{(D,w)}^+(v) \le f(v)$ for each vertex $v$. The arc-weighted   degeneracy $d_{w}(G)$ of $G$ is the minimum $d$ such that $G$ is  arc-weighted   $d$-degenerate.  
   Given an arc-weighted acyclic orientation $(D,w)$ of $G$ with $d_{(D,w)}^+(v) \le  f(v)$,  there is not only an easy algorithm that constructs an $(L,M)$-colouring of $G$ for any $(f+1)$-DP-cover $(L,M)$ of $G$, but also an easy winning strategy for the DP-$(f+1)$-painting game on $G$. Moreover, arc-weighted  $f$-degeneracy implies $(f+1)$-Alon--Tarsi. As an arc-weighted acyclic orientation $(D,w)$ of $G$ with $w(e)=1$ for all $e$ is equivalent to a (non-weighted)  acyclic orientation of $G$, $d_{w}(G)$ is bounded from above by the degeneracy $d(G)$ of $G$. We observe that the difference $d(G)-d_{w}(G)$ can be arbitrarily large.
 Thus, $d_{w}(G)+1$ provides a better upper bound than $d(G)+1$ on its DP-paint number (and hence its DP-chromatic number, paint number and choice number), as well as its Alon--Tarsi number. Thomassen's proof of the 5-choosability of planar graphs can be easily adapted to prove  that all planar graphs are arc-weighted   $4$-degenerate, implying that planar graphs are DP-$5$-paintable as well as $5$-Alon--Tarsi. As an extension of the characterization of degree-choosable graphs, we prove that a connected graph $G$ is arc-weighted  $(d_G-1)$-degenerate unless $G$ is a GDP-tree. In particular, $d_{w}(G) \le \Delta(G)-1$ unless $G$ is a complete graph or a cycle.  
   Then we prove that 3-connected non-complete planar graphs are degree-truncated arc-weighted   $14$-degenerate, implying that such graphs have degree-truncated DP-paint number at most $15$. The previous known upper bound for this parameter was $16$.   
	\end{abstract}
	
	\section{Introduction}	
	
	For a graph $G$, an orientation $D$ of $G$ is a digraph  obtained from $G$ by replacing each edge $uv$ with an arc $(u,v)$ or $(v,u)$.    A digraph  $D$
is   {\em acyclic} if it contains no directed cycles. The acyclic orientation of  graphs is a fundamental concept in graph theory and has numerous applications. 

	Given a digraph $D$ and a vertex $v$ of $D$, $N_D^+(v), N_D^-(v)$ are the sets of out-neighbours and in-neighbours of $v$, respectively. The {\em out-degree} of $v$ is $d_D^+(v)=|N_D^+(v)|$, and the {\em in-degree} of $v$ is $d_D^-(v)=|N_D^-(v)|$. 

We denote by $\mathbb{N}^G$ the set of mappings $f: V(G) \to \mathbb{N}=\{0,1,\ldots\}$. For $f,g \in \mathbb{N}^G$ and $a, b \in \mathbb{N}$,  $af+ bg \in \mathbb{N}^G$  is defined as $(af+bg)(v)=af(v)+bg(v)$.  We write $f \le g$ if $f(v) \le g(v)$ for all $v \in V(G)$. For brevity, for an integer $k$, we denote by $k$ the constant mapping in $\mathbb{N}^G$ defined as $k(v)=k$ for all vertices $v$.  We denote by $d_G(v)$ the degree of $v$, and use $d_G$ as a mapping in $ \mathbb{N}^G$.   

	\begin{definition}
		\label{def-degenerate}
		For $f \in \mathbb{N}^G$, we say $G$ is \emph{$f$-degenerate} if there is an acyclic orientation $D$ of $G$ such that 
		$d_D^+(v) \le  f(v)$ for each vertex $v$.     
  The \emph{degeneracy} $d(G)$ is the minimum integer $k$ such that $G$ is  $k$-degenerate.  
  \end{definition}

	The number  $d(G)+1$ is called the  {\em colouring number} of $G$ (also called the {\em strict degeneracy} of $G$), and  is an upper bound for many colouring parameters. We shall review the definitions of those parameters in Section~2. Here, for an illustration of these relations,  we  define   list colouring and choice number   of a graph. 
 
 \begin{definition}
    \label{def-list}
    Let $G$ be a graph. A {\em list assignment} of $G$ is a mapping $L$ which assigns to each vertex $v$ of $G$ a set $L(v)$ of permissible colours. Given  a list assignment $L$ of $G$, an {\em $L$-colouring} of $G$ is a mapping $\phi: V(G) \to \cup_{v \in V(G)}L(v)$ such that for each vertex $v$, $\phi(v) \in L(v)$ and for any edge $e=uv$ of $G$, $\phi(u)\ne \phi(v)$. For $f \in \mathbb{N}^G$, an {\em $f$-list assignment} of $G$ is a list assignment $L$ of $G$ with $|L(v)| \ge f(v)$ for each vertex $v$. We say $G$ is {\em $f$-choosable} if $G$ is $L$-colourable for every $f$-list assignment $L$. The {\em choice number} $\operatorname{ch}(G)$ of $G$ is the minimum $k$ for which $G$ is $k$-choosable.
\end{definition}
 
  Assume $G$ is $f$-degenerate,  $D$ is an acyclic orientation  of  $G$ with $d_D^+(v ) \le f(v )$ and  $L$ is an $(f+1)$-list assignment of $G$. Then, $D$ induces an ordering $v_1,v_2, \ldots,v_n$ of the vertices of $G$ such that  an edge $v_iv_j$ of $G$ is oriented as $(v_i, v_j)$ if and only if $j < i$. Since $v_i$ has only $d_D^+(v_i )$  neighbours $v_j$ with $j < i$  and $|L(v_i)| \ge f(v_i)+1 > d_D^+(v_i)$, we can properly colour the vertices of $G$ greedily in the order $v_1,v_2, \ldots, v_n$  with colours from their lists. Therefore $G$ is $(f+1)$-choosable.

\begin{definition}
    \label{def-arcweighted}
    An {\em arc-weighted orientation}  of a graph $G$ is a pair $(D,w)$, where $D$ is an orientation of $G$   and $w: E(D) \to \{1,2,\ldots\}$ is a mapping that assigns to each arc of $D$ a positive integer as its weight.
\end{definition}
 
Given an arc-weighted orientation $(D,w)$ of a graph $G$, for each vertex $v$, the {\em weighted out-degree} and {\em weighted in-degree} of $v$ are defined as
$$d_{(D,w)}^+(v) = \sum_{u \in N_D^+(v)}w((v,u)), \quad d_{(D,w)}^-(v) = \sum_{u \in N_D^-(v)}w((u,v)).$$

A natural problem is how to generalize the concept of acyclic orientation to arc-weighted orientations, and provide better upper bounds for various colouring parameters. 

An orientation $D$ of $G$ can be viewed as an arc-weighted orientation $(D,w)$ with $w(e)=1$ for each arc $e$. A digraph $D$ is  acyclic if and only if every sub-digraph $D'$ has a sink, i.e., a vertex $v$ with $d_{D'}^+(v)=0$. If $D'$ has at least one arc, then one of the arcs $e=(u,v)$ has its head $v$ being a sink, and hence $w(e) > d_{D'}^+(v)$. This motivates the following definition of a dominating arc in an arc-weighted digraph.

\begin{definition}
    \label{def-dominating}
    Assume $(D,w)$ is an arc-weighted digraph. An arc $e=(u,v)$ is called {\em dominating} if $w(e) > d_{(D,w)}^+(v)$. 
\end{definition}

An arc $e=(u,v)$ with $w(e)=1$ is dominating if and only if $v$ is a sink. Thus one can alternatively  define acyclic digraphs (viewed as arc-weighted digraphs with each arc of weight $1$) as follows: 
A digraph $D$ is acyclic if every non-empty sub-digraph $D'$ of $D$ has a dominating arc. This leads to the following generalization of acyclic digraphs to arc-weighted acyclic digraphs.

\begin{definition}
    \label{def-strong}
    An arc-weighted digraph $(D,w)$ is {\em acyclic} if, for every non-empty sub-digraph $D'$ of $D$, $(D', w)$ contains a dominating arc. We say that a graph $G$ is {\em arc-weighted $f$-degenerate} if $G$ has an arc-weighted acyclic orientation $(D,w)$ with $d_{(D,w)}^+(v) \le f(v)$ for each vertex $v$.
\end{definition}

We interpret a directed cycle $C$ in a digraph $D$ as  a "base" for "circulations of flow" in $D$. In this sense, a directed cycle $C$ in $D$ is not a directed cycle in $(D,w)$ if it contains a dominating arc $e=(u,v)$: the flow of weight $w(e)$ is "blocked" at the vertex $v$. So we use the term acyclic for such an arc-weighted digraph $(D,w)$. Some other generalizations of acyclic orientation to arc-weighted orientation are discussed in Section 4.

\begin{definition}
    The {\em arc-weighted   degeneracy} $d_{w}(G)$ of $G$ is the minimum integer $k$ such that $G$ is arc-weighted  $k$-degenerate.  
\end{definition}

 Assume $|E(D)|=m$. 
Let $D_0=D$, and for $i=1,2,\ldots, m$, let $e_i=(x_i,y_i)$ be a dominating arc in
$D_{i-1}$ and $D_i=D_{i-1} - e_i$. We obtain an ordering  $e_1=(x_1,y_1), e_2=(x_2,y_2), \ldots, e_m=(x_m,y_m)$ of the arcs of $D$  such that 
$$w(e_i) > d_{(D-\{e_1,e_2,\ldots, e_{i-1}\}, w)}^+(y_i) \eqno{(1)}.$$
Conversely, if the arcs of $D$ can be ordered this way, then every non-empty sub-digraph $D'$ of $D$ has a dominating arc (with respect to $w$). So the arcs of $D$ having an ordering satisfying (1) is an equivalent definition of $(D,w)$ being acyclic.

We say $(D,w)$ is {\em minimal acyclic} if $w(e_i) = d_{(D-\{e_1,e_2,\ldots, e_{i-1}\}, w)}^+(y_i)+1$ for $i=1,2,\ldots, m$ in (1).
Note that if $(D,w)$ is acyclic, then there is a weighting function $w'$ such that $w' \le w$ and $(D,w')$ is minimal acyclic. For the purpose of deriving upper bounds for colouring parameters, it suffices to consider minimal arc-weighted acyclic digraphs.

\begin{figure}
    \centering
    \includegraphics[width=0.26\textwidth]{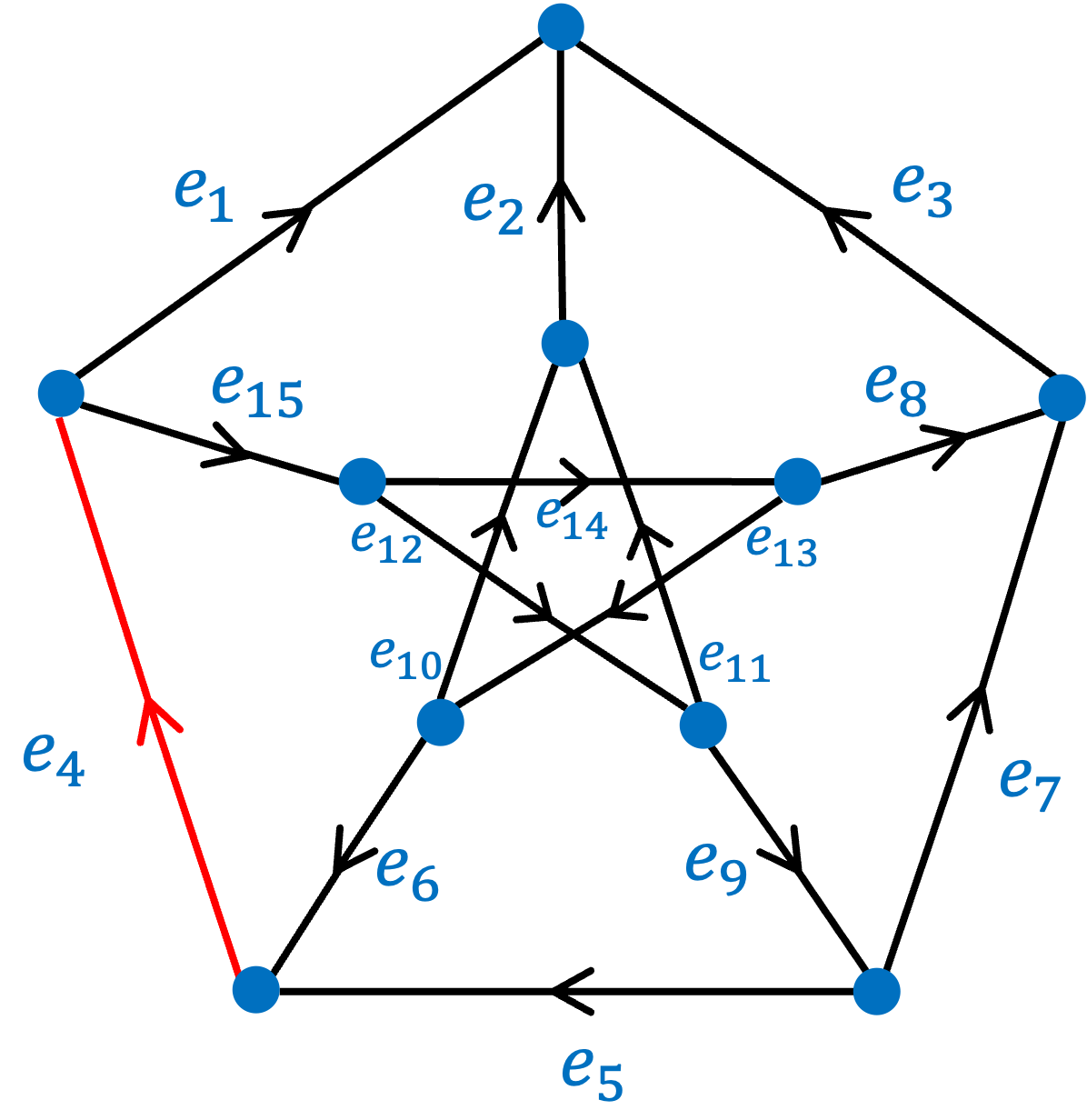}
    \caption{An arc-weighted acyclic orientation of the Petersen graph, where $w(e)=1$ for a black arc and $w(e)=2$ for a red arc. For each $i$, $e_i$ is a dominating arc in $D-\{e_1,\ldots, e_{i-1}\}$.}
    \label{fig-petersen}
\end{figure}

It follows from the definition that $d_w(G)\le d(G)$ for any graph $G$.
Figure \ref{fig-petersen} shows that the Petersen graph $P$ has $d_{w}(P)\le 2$, while $d(P)=3$. We shall observe later that the difference $d(G) - d_{w}(G)$ can be arbitrarily large.

The degeneracy $d(G)$ gives an upper bound $d(G)+1$ for many colouring parameters of $G$, including the DP-paint number and the Alon--Tarsi  number of $G$. We shall show that $d_{w}(G)+1$ is also an upper bound for the DP-paint number and the Alon--Tarsi  number of $G$. In the following, we prove that the arc-weighted $f$-degeneracy implies the $(f+1)$-choosability, and  the proof illustrates the idea of finding a proper colouring  via a given  arc-weighted acyclic orientation $(D,w)$.

     For  a subset $X$ of $V(G)$, $1_X$ is the characteristic function of $X$, defined as 
  \[
  1_X(v)=\begin{cases}
      1, &\text{ if $v \in X$}, \cr 
      0, &\text{ otherwise.}
  \end{cases}
  \]
  For simplicity, we write $1_x$ for $1_{\{x\}}$. 
  
\begin{lem}
    \label{lem-choosable}
    If $G$ is arc-weighted   $f$-degenerate, then $G$ is $(f+1)$-choosable.
    \end{lem}
    \begin{proof}
       Assume $G$ is arc-weighted   $f$-degenerate, and $L$ is an $(f+1)$-list assignment of $G$. We prove  that $G$ is $L$-colourable. The proof is based on the induction of the number of edges of $G$. 
       
If $E(G)=\emptyset$, then choosing any colour $\phi(v)\in L(v)$ for each vertex $v$ gives an $L$-colouring of $G$.
 
  Assume $|E(G)| \ge 1$, and $(D,w)$ is a  minimal arc-weighted acyclic orientation of $G$ with $d_{(D,w)}^+   \le  f$. Without loss of generality, we may assume that $|L(v)|=f(v)+1=d_{(D,w)}^+(v)+1 $  for each vertex $v$. 

By definition, $(D,w)$ has a dominating arc $e=(x,y)$.  Let $G'=G-xy$, let $f'=   f-(f(y)+1)1_x$, and let 
  $L'$ be the list assignment of $G'$ obtained from $L$ by deleting   $L(y)$ from $L(x)$, i.e., $L'(x)=L(x)-L(y)$ and $L'(v)=L(v)$ for $v \ne x$. Since $|L(y)|=f(y)+1$, we have 
  $|L'(x)| \ge |L(x)|-|L(y)| =f(x)+1-(f(y)+1)=f'(x)+1$. Hence $L'$ is an $(f'+1)$-list assignment of $G'$. 
  
  Let $D'=D-\{(x,y)\}$, and let $w'$ be the restriction of $w$ to $E(D')$. Since $e$ is dominating and $f(y)=d_{(D,w)}^+(y)$, we have $w(e)\ge f(y)+1$. Therefore,
$$d_{(D',w')}^+(x)=d_{(D,w)}^+(x)-w(e)
\le f(x)-(f(y)+1)=f'(x),$$
and for every $v\ne x$,
$d_{(D',w')}^+(v) \le f(v)=f'(v)$.
Thus $d_{(D',w')}^+\le f'$.
 Moreover, $(D',w')$ is arc-weighted acyclic. According to the induction hypothesis, $G'$ has a $L'$-colouring $\phi$, which is also an $L$-colouring of $G$.   
    \end{proof}

The proof of Lemma \ref{lem-choosable} gives an algorithm that finds an $L$-colouring of $G$ via an arc-weighted acyclic orientation $(D,w)$ with $d_{(D,w)}^+(v)< |L(v)|$ for each vertex $v$. The colouring procedure is to remove edges one by one, to obtain an empty graph $G'$ on the same vertex set,
and with each vertex having at least one colour available. Moreover, a proper colouring of $G'$ using the available colours is also a proper colouring of $G$. It is interesting to compare this colouring process with the greedy colouring process that is used to derive a proper colouring through an acyclic orientation of $G$. 

Let $\mathcal{L}$ be the set of pairs $(G,f)$, where $G$ is a graph, and $f \in \mathbb{N}^G$. For convenience, we let $(\emptyset, \emptyset)$ be a pair in $\mathcal{L}$, which stands for the non-vertex graph (the term "empty graph" is reserved for an edgeless graph). For $f \in \mathbb{Z}^G$ and a subset $U$ of $V(G)$, let $f|U$ be the restriction of $f$ to $U$, and let  $f_{-U}: V(G)-U \to \mathbb{Z}$ be defined as 
		$f_{-U}(x) = f(x) - |N_G(x) \cap U|$ for $x \in V(G)-U$ and $f_{-U}=\emptyset$ if $U=V(G)$.
We write $f_{-v}$ for $f_{-\{v\}}$, and 
  when there is no confusion, we may use   $f$ for $f|U$.
		For any $u\in V(G)$,    $E_G(u)$  is the set of edges incident to  $u$ in $G$.

Given a pair $(G,f) \in \mathcal{L}$ and a vertex $x$ of $G$, the {\em vertex-deletion  operation} $\dd_x$ on $(G,f)$ is defined as   
	$$\dd_x(G,f)=(G-x,f_{-x}). $$
We say that operation $\dd_x$ is {\em legal} for $(G,f)$ if $f_{-x} (v) \ge 0$ for   $v \in V(G-x)$. It is easy to see that a graph $G$ is   $f$-degenerate if, starting from $(G,f)$, there is a sequence of legal vertex-deletion operations that delete all vertices of $G$. We may view the move of deleting a vertex $x$ as colouring the vertex $x$: Originally, each vertex $x$ has $f(x)+1$ permissible colours. Once $x$ is coloured, each of its neighbors loses one permissible colour. The legality condition ensures that no remaining vertex loses more colours than it can afford.

In the proof of Lemma \ref{lem-choosable}, in removing the dominating arc $e=(x,y)$,
we remove all the colours in $L(y)$ from $L(x)$. Formally, we define
the operation of arc-deletion as follows:

\begin{definition}
    \label{def-ed}
     Given $(G, f) \in \mathcal{L}$, if $xy \in E(G)$  and $s \ge 0$ is an integer, the {\em arc-deletion} operation  $\ed_{(x,y),s}$ is defined as $$\ed_{(x,y),s}(G,f)=(G', f'),$$
     where $G'=G-xy$ and 
     \[
     f'(v) = \begin{cases}
         f(y)-s, &\text{ if $v=y$} \cr 
         f(x) - (f(y)+1-s), &\text{ if $v=x$}, \cr 
         f(v), &\text{ otherwise.}
     \end{cases}
     \]
     The operation  $\ed_{(x,y),s}$ is legal for $(G,f)$ if $f(x) \ge  f(y)+1-s$ and $f(y) \ge s$. 
\end{definition}

 For brevity, we write $\ed_{(x,y)}$ for $\ed_{(x,y), 0}$. 

 The motivation for this definition is as follows: Assume $L$ is an $(f+1)$-list assignment of $G$. To delete the edge $xy$, we modify the lists at its two ends so that they become disjoint. First we delete $s$ colours from $L(y)$ to obtain the list $L'(y)$, and then delete all colours in $L'(y)$ from $L(x)$ to obtain the list $L'(x)$. For $v \ne x,y$, let $L'(v)=L(v)$.  
Since $L'(x)\cap L'(y)=\emptyset$, every $L'$-colouring of $G-xy$ is an $L$-colouring of $G$. The legality conditions mean precisely that the two endpoints can afford these losses: $y$ loses $s$ colours, while $x$ is allowed to lose up to $f(y)+1-s$ colours.
Thus if $G-xy$ is $(f'+1)$-choosable, then $G$ is $(f+1)$-choosable.

 We may recursively apply a sequence of arc-deletion operations to a pair $(G,f) \in \mathcal{L}$. Let
   $$ \Omega=\ed_{(x_1,y_1),s_1}\cdots \ed_{(x_k,y_k),s_k}$$
be a sequence of arc-deletion operations. If $k=0$, then $\Omega$ is the
identity operation. If $k\ge 1$, let
    $$\Omega'=\ed_{(x_1,y_1),s_1}\cdots
    \ed_{(x_{k-1},y_{k-1}),s_{k-1}}.$$
We define
$$ \Omega(G,f)=\Omega'\bigl(\ed_{(x_k,y_k),s_k}(G,f)\bigr).$$

We say that $\Omega$ is legal for $(G,f)$ if
$\ed_{(x_k,y_k),s_k}$ is legal for $(G,f)$, and $\Omega'$ is legal for $\ed_{(x_k,y_k),s_k}(G,f)$.

The following lemma gives an alternate definition of arc-weighted $f$-degeneracy. Given a graph $G$, we denote by $G_0$ the edgeless graph with $V(G_0)=V(G)$.  

\begin{lem}
    \label{lem-alternate}
    A graph $G$ is arc-weighted $f$-degenerate if and only if there is a sequence $\Omega$ of legal arc-deletion operations such that $\Omega(G,f)=(G_0,f_0)$ and $f_0 \ge 0$.  
\end{lem}
\begin{proof}
     Assume $G$ is arc-weighted $f$-degenerate and $(D,w)$ is an arc-weighted  acyclic orientation of $G$ with $d_{(D,w)}^+  \le f$. 
     We prove by induction on the number of edges of $G$ that  there is a sequence $\Omega$ of legal arc-deletion operations such that $\Omega(G,f)=(G_0,f_0)$ and $f_0 \ge 0$.   
     If $G$ is edgeless, then there is nothing to prove. Otherwise, there is an arc $e=(u,v)$ such that $w(e) \ge d_{(D,w)}^+(v)+1$.  
     
        Let $s = f(v) - d_{(D,w)}^+(v)$, and let $(G-e,f')=\ed_{(u,v),s}(G,f)$.  
      It follows from the definition that $G-e$ is  arc-weighted $f'$-degenerate. 
     By the induction hypothesis, there is a sequence $\Omega'$ of arc-deletion operations  such that $\Omega'(G-e, f') = (G_0,f_0)$ with $f_0 \ge 0$. Let $\Omega = \Omega' \ed_{(u,v),s}$.   Then $\Omega$ is a sequence of legal moves such that $\Omega(G, f) = (G_0,f_0)$ with $f_0 \ge 0$.

 Conversely, assume that there is a sequence $\Omega$ of legal moves such that $\Omega(G,f) = (G_0,f_0)$ with $f_0 \ge 0$. 
     We show by induction on $\sum_{v \in V(G)}f(v)$ that there is an arc-weighted  acyclic orientation $(D,w)$ of $G$ with 
    $d_{(D,w)}^+(v) \le f(v)$ for each vertex $v$. Again, if $G$ is edgeless,   then there is nothing to prove. Assume that the first move is  $\ed_{(u,v), s}$, $\ed_{(u,v), s}(G, f) = (G-uv, f')$, and $\Omega (G, f)= \Omega' \ed_{(u,v), s} (G, f) = \Omega'(G-uv, f')$.
 By the induction hypothesis, 
    there is an arc-weighted acyclic orientation $(D',w')$ of $G-uv$ such that 
    $d_{(D',w')}^+(x) \le f'(x)$ for all $x \in V(G)$. We extend $(D',w')$ to an arc-weighted orientation $(D,w)$ by adding the arc $e=(u,v)$ with weight $w(e) = f'(v)+1$. Then $(D,w)$ is  arc-weighted acyclic  and $d^+_{(D,w)}(x)   \le f(x)$ for each vertex $x$. 
\end{proof}

The vertex-deletion operation has the same effect as a  sequence of arc-deletion operations: Assume $N_G(x)=\{y_1,y_2, \ldots, y_k\}$. Let $s= f(x)$.  Then $\dd_x$ has the same effect as 
$$\ed_{(y_k,x)}\ed_{(y_{k-1},x)} \ldots \ed_{(y_2,x)} \ed_{(y_1,x), s}.$$ 
Note that after this sequence of operations, $x$ becomes an isolated vertex, which can be viewed as deleted.
On the other hand, the arc-deletion operation is more flexible and will lead to a better upper bound for many colouring parameters. 

In Section 2, we define some colouring parameters, and prove that the arc-weighted degeneracy plus one provides an upper bound for all these colouring parameters. 

In Section 3, we adapt Thomassen's proof to show that every planar graph is arc-weighted $4$-degenerate, and hence $5$-DP-paintable and $5$-AT.

In Section 4, we study the degree-truncated arc-weighted strict degeneracy of a graph. We prove that every 3-connected non-complete planar graph has degree-truncated strict degeneracy at most $15$. This implies that every 3-connected non-complete planar graph is degree-truncated $15$-DP-colourable, and degree-truncated $15$-AT, improving the previous known upper bounds on these parameters by $1$. 

In Section 5, we make some remarks and pose some open problems. In particular, we compare arc-weighted degeneracy with weak degeneracy introduced by Bernshteyn and Lee \cite{BL}, show that the difference $d(G) - d_{w}(G) $ can be arbitrarily large.

	\section{ $f$-Alon--Tarsi and  $f$-DP-paintable  graphs} 
    \label{f-AT-and-f-DP-paintable}
	For a graph $G$, $d(G)+1$  is an upper bound for many colouring parameters, including the DP-paint number and the Alon--Tarsi number of $G$. We shall show that   $d_{w}(G)+1$ is also an upper bound for its DP-paint number and Alon--Tarsi number.

	 In this section, we first review the definitions of some colouring parameters. Then we prove that if a graph $G$ is arc-weighted $f$-degenerate, then $G$ is $(f+1)$-Alon--Tarsi and  $(f+1)$-DP-paintable.

For a graph $G$, let $<$ be an arbitrary ordering of $V(G)$ and let $$P_G(x_v: v \in V(G)) = \prod_{uv \in E(G), u< v}(x_u-x_v).$$
For $t \in \mathbb{N}^G$, let $c_{G,t}$ be the coefficient of the monomial $\prod_{v \in V(G)}x_v^{t(v)}$ in the expansion of $P_G$. For $f \in \mathbb{N}^G$, we say that $G$ is {\em $f$-Alon--Tarsi} ($f$-AT for short) if for some $t \in \mathbb{N}^G$ with $t+1  \le f$, $c_{G,t} \ne 0$. The {\em Alon--Tarsi number} of $G$ is defined as 
$$\operatorname{AT}(G)=\min\{k: G \text{ is $k$-AT}\}.$$

For a mapping $\phi: V(G) \to \mathbb{N}$, let $P_G(\phi)$ be the evaluation of $P_G$ at $x_v = \phi(v)$ for $v \in V(G)$. Then $P_G(\phi) \ne 0$ if and only if $\phi$ is a proper colouring of $G$, i.e., $\phi(u) \ne \phi(v)$ for each edge $uv$ of $G$. Given a list assignment $L$ of $G$, finding a proper $L$-colouring of $G$ is equivalent to finding a non-zero point of $P_G$ in the grid $ \bigtimes_{v \in V(G)} L(v)$. Combinatorial Nullstellensatz \cite{AT} gives a sufficient condition for the existence of a non-zero point of a polynomial $f(x_1,x_2,\ldots,x_n)$ in a grid $S_1 \times S_2 \times \ldots \times S_n$: If $\prod_{i=1}^n x_i^{t_i}$ is a highest degree monomial of $f$,   the coefficient   of $\prod_{i=1}^n x_i^{t_i}$
in the expansion of $f$ is non-zero and $t_i \le |S_i|-1$, then $f$ has a non-zero point   in   $S_1 \times S_2 \times \ldots \times S_n$. From this it follows that $\operatorname{ch}(G) \le \operatorname{AT}(G)$ for every graph $G$. 

Recently, Kozik and Podkanowicz considered a polynomial associated with edge-weighted graphs and gave an alternate definition of Alon--Tarsi number \cite{kozik}.

Let $w \in \mathbb{N}^{E(G)}$ be an integer weighting of the edges of $G$ such that $w(e) \ge 1$ for each edge $e$. Let $$P_{G,w}(x_v: v \in V(G)) = \prod_{e=uv \in E, u< v}(x_u^{w(e)} - x_v^{w(e)}).$$ 
For $t \in \mathbb{N}^G$, let $c_{(G,w),t}$ be the coefficient of the monomial $\prod_{v \in V(G)}x_v^{t(v)}$ in the expansion of $P_{G,w}$. 

It is easy to see that $P_{G,w} = P_G \cdot g$ for some polynomial $g$. Hence, if $c_{(G,w),t} \ne 0$, then there exists $t' \in \mathbb{N}^G$ such that $t' \le t$ and $c_{G,t'} \ne 0$. Thus we have the following lemma.

\begin{lem}
    \label{lem-AT}
    A graph $G$ is $f$-AT if and only if there exist a weighting $w \in \mathbb{N}^{E(G)}$  of the edges of $G$ so that $w(e) \ge 1$ for each edge $e$, and $t \in \mathbb{N}^G$ such that $c_{(G,w),t} \ne 0$ and $t+1 \le f$. 
\end{lem}  

\begin{definition}
     An arc-weighted digraph $(D,w)$ is called {\em Eulerian} if $d^+_{(D,w)}(v) = d_{(D,w)}^-(v)$ for
 each vertex $v$.
\end{definition}
Assume $(D,w)$ is an arc-weighted orientation  of $G$. Let 
\begin{eqnarray*}
    \mathcal{E}(D,w) &=& \{H \subseteq E(D): (D[H],w) \text{ is an Eulerian weighted digraph}\},\\
    \mathcal{E}_e(D,w) &=& \{H \in  \mathcal{E}(D,w): |H| \text{ is even}\},\\
     \mathcal{E}_o(D,w) &=& \{H \in  \mathcal{E}(D,w): |H| \text{ is odd}\},\\
     {\rm diff}(D,w) &=& |\mathcal{E}_e(D,w)|-|\mathcal{E}_o(D,w)|.
\end{eqnarray*}
 It is known \cite{kozik} that 
$$c_{(G,w), d_{(D,w)}^+} = \pm {\rm diff}(D,w).$$
We say $(D,w)$ is an {\em arc-weighted AT-orientation} of $(G,w)$ if ${\rm diff}(D,w) \ne 0$. If $(D,w)$ is an   arc-weighted AT-orientation   of $(G,w)$ and $w(e)=1$ for each edge $e$, then we say $D$ is an AT-orientation of $G$.
Thus, we have the following lemma.
\begin{lem}
    \label{lem-3}
 For a graph $G$ and $f \in \mathbb{N}^G$, the following are equivalent:
    \begin{enumerate}
        \item $G$ is $f$-AT.
        \item $G$ has an   AT-orientation $D$ with   $d_D^+ +1 \le f$.  
        \item $G$ has an arc-weighted AT-orientation $(D,w)$  with $d_{(D,w)}^+ +1 \le f$.
    \end{enumerate}
\end{lem}

\begin{theorem}
    \label{thm-AT}
    If $G$ has an arc-weighted  acyclic orientation $(D,w)$, then $G$ is $(d_{(D,w)}^++1)$-AT. So if $G$ is arc-weighted   $f$-degenerate, then $G$ is $(f+1)$-AT.
\end{theorem}
\begin{proof}
For any non-empty sub-digraph $D'$ of $D$, $(D', w)$ has a dominating arc. Hence $(D',w)$ is not an arc-weighted  Eulerian sub-digraph.  Thus, $(D,w)$ has only one arc-weighted Eulerian sub-digraph, namely, the empty sub-digraph. 
Thus, ${\rm diff}(D,w)\neq 0$, and  
   $(D,w)$ is an arc-weighted AT-orientation of $G$. 
\end{proof}

Next, we show that if $G$ is arc-weighted  $f$-degenerate, then $G$ is $(f+1)$-DP-paintable. First we review the definition of  DP-colouring, painting game, and DP-painting game.

 \begin{definition}
	\label{def-cover}
	A {\em{\rm DP}-cover} of a graph $G $ is a pair $(L,M)$, where $L$ is a mapping that assigns to each vertex $v$ a set $L(v)$ of permissible colours such that  $L(v) \cap L(u) = \emptyset$ for $u \ne v$, and $M=\{M_{e}: e \in E(G)\}$, where for each edge $e=uv$, $M_e$ is a   matching between $L(u)$ and $L(v)$. For  $f \in \mathbb{N}^G$, we say $(L,M)$ is an $f$-DP-cover of $G$ if $|L(v)|\ge f(v)$ for each vertex $v \in V(G)$.
\end{definition}

\begin{definition}
	\label{def-coloring}
	Given a DP-cover $(L,M)$ of a graph $G$, an $(L,M)$-colouring of $G$ is a mapping $\phi: V(G) \to \bigcup_{v \in V(G)}L(v)$ such that for each vertex $v \in V(G)$, $\phi(v) \in L(v)$, and for each edge $e=uv \in E(G)$, $\phi(u)\phi(v) \notin  M_e$. We say $G$ is {\em $(L, M)$-colourable} if it has an $(L,M)$-colouring.
\end{definition}

 For an $(L,M)$-colouring $\phi$ of $G$, we say the vertex $v$ of $G$ is coloured by the colour $\phi(v)$.

\begin{definition}
	\label{DP-coloring}
	Let $G$ be a graph and $f \in \mathbb{N}^{G}$.  We say $G$ is {\em $f$-DP-colourable} if for every $f$-DP-cover $(L,M)$, $G$ has an $(L,M)$-colouring. The {\em DP-chromatic number} of $G$ is defined as 
 $$\chi_{\mathrm{DP}}(G)=\min\{k: G \text{ is $k$-DP-colourable} \}.$$
\end{definition}

Given an $f$-list assignment $L$ of a graph $G$, let $(L',M)$ be the  $f$-DP-cover of $G$ induced by $L$, where  $L'(v) = \{(i,v): i \in L(v)\}$ for each vertex $v$ of $G$, and $M=\{M_{uv}: uv \in E(G)\}$ is defined as $$M_{uv}=\{\{(i,u)(i,v)\}: i \in L(u) \cap L(v)\}$$
for each edge $uv$ of $G$.
It is obvious that $G$ is $L$-colourable if and only if $G$ is $(L',M)$-colourable. Therefore  if $G$ is $f$-DP-colourable, then it is $f$-choosable, and hence $\operatorname{ch}(G) \le \chi_{\mathrm{DP}}(G)$. The proof of Lemma \ref{lem-choosable} can be easily adapted to prove the following result:

\begin{lem}
    \label{lem-dpf}
    If $G$ is arc-weighted $f$-degenerate, then $G$ is  $(f+1)$-DP-colourable. 
\end{lem}

The DP-paint number $\chi_{\mathrm{DPP}}(G)$ is the online version of DP-chromatic number, and is defined through a two-player game  \cite{KKLZ2020}.

 \begin{definition}
	\label{def-coloring2}
	Given a DP-cover $(L,M)$ of a graph $G$, a {\em partial $(L,M)$-colouring} of $G$ is a mapping $\phi: W \to \bigcup_{v \in V(G)}L(v)$ for a subset $W$ of $V(G)$ such that for each vertex $v \in W$, $\phi(v) \in L(v)$, and for each edge $e=uv \in E(G[W])$, $\phi(u)\phi(v) \notin  M_e$. We say $W$ is the set of coloured vertices in this partial $(L,M)$-colouring.  
\end{definition}

\begin{definition}
    \label{def-dppainting}
	For $f \in \mathbb{N}^G$, the \emph{$f$-DP-painting game} on $G$ is played by two players: Lister and Painter. Initially each vertex $v$ has $f(v)$ tokens and is uncoloured. 
  
  In the $i$th round, Lister chooses a cover $(L_i,M_i)$ of $G$, and  for each vertex $v$,  removes $|L_i(v)|$ tokens from $v$ (thus $|L_i(v)|$ is at most the number of remaining tokens of $v$).   Painter constructs a partial $(L_i, M_i)$-colouring $h_i$ of $G$ that colours a subset $X_i$ of vertices of $G$.  The game ends either if all vertices of $G$ are coloured, in which case Painter wins the game, or some uncoloured vertex has no token left, in which case Lister wins. We say $G$ is \emph{$f$-DP-paintable} if Painter has a winning strategy for this game.  
		The \emph{DP-paint number}	$\chi_{\mathrm{DPP}}(G)$ of $G$ is the minimum $k$ such that $G$ is $k$-DP-paintable.
\end{definition}	

If Lister chooses an $f$-DP-cover $(L,M)$ in the first move, then Painter needs to respond with an $(L,M)$-colouring of $G$ (instead of a partial $(L,M)$-colouring of $G$), for otherwise, Painter loses the game. So if $G$ is $f$-DP-paintable, then $G$ is $f$-DP-colourable, and the converse is not true. 
		
In the definition of $f$-DP-painting game, if each cover $(L_i,M_i)$ chosen by Lister is required to satisfy the additional constraint that $|L_i(v)| \le 1$ for each vertex $v$, then the game is called the {\em $f$-painting game}.  This definition of the $f$-painting game is slightly different from (but equivalent to) the standard formulation in the literature (see, e.g., \cite{Schauz}). We say $G$ is {\em $f$-paintable} if Painter has a winning strategy in the $f$-painting game on $G$, and the {\em paint number} $\chi_P(G)$ is the minimum $k$ such that $G$ is $k$-paintable.
Equivalently, in the $i$th round, Lister reveals a set $U_i = \{v: |L_i(v)|=1\}$ of uncoloured vertices, which can be understood to be the set of vertices $v$ whose lists contain colour $i$. Painter chooses an independent set $X_i$ of $U_i$, and colours vertices of $X_i$ by the colour $i$. In this sense, the $f$-painting game on $G$ is an online version of $f$-list colouring of $G$.

  It follows from the definitions that for every graph $G$, $\operatorname{ch}(G) \le \chi_P(G) \le \chi_{\mathrm{DPP}}(G)$ and 
  $\chi_{\mathrm{DP}}(G) \le \chi_{\mathrm{DPP}}(G)$. 
		
The following lemma in \cite{KKLZ2020} can be viewed as an alternate definition of $f$-DP-paintability and $f$-paintability.

 \begin{lem}
 \label{lem-0}
    Let $G$ be a graph and $f \in \mathbb{N}^{G}$.  Then $G$ is $f$-DP-paintable (respectively, $f$-paintable) if and only if one of the following holds:
  \begin{itemize}
  \item[(i)] $V(G)=\emptyset$.
  \item[(ii)] For every $g$-DP-cover $(L,M)$ of $G$ with $0 \ne g \le f$ (respectively, for every $g$-DP-cover $(L,M)$ of $G$ with $g(v) \le \min\{1, f(v)\}$ for each vertex $v$), there is a subset $X$ of $V(G)$, so that $G[X]$ is $(L,M)$-colourable, and $G-X$ is $(f-g)$-DP-paintable (respectively, $(f-g)$-paintable).
  \end{itemize} 
 \end{lem}


 Note that for a graph $G$ to be $f$-DP-colourable, each vertex $v$ has at least 1 colour, i.e., $f(v) \ge 1$. In the DP-painting game, we consider DP-colouring of subgraphs of $G$. The mappings $g,f$ involved in the description of the game may take value $0$ on some vertices.

		\begin{lem}
			\label{lem-1}
			Let $G$ be a graph, $f,g \in \mathbb{N}^G$ and $0 \ne g \le f+1$. If $G$ is arc-weighted  $f$-degenerate, then  there is a non-empty subset $X$ of $V(G)$ such that $G[X]$ is arc-weighted  $(g-1)$-degenerate and $G-X$ is arc-weighted $(f-g)$-degenerate. 
		\end{lem}
		\begin{proof}
	The proof is by induction on $\sum_{v \in V(G)}f(v)$. 

 
 Assume that $G$ is arc-weighted $f$-degenerate. 
 If $E(G)=\emptyset$, then let $X=\{v: g(v)  \ge 1\}$, and the conclusion holds.  
 
  Assume $|E(G)| \ge 1$, and $(D,w)$ is a  minimal arc-weighted   acyclic orientation of $G$. If there is a vertex $x$ for which  
$f(x) > d_{(D,w)}^+(x)$, 
then let $f'=f-1_{\{x\}}$ and $g' = \max\{0, g- 1_{\{x\}}\}$. Then $g' \le f'+1$ and $G$ is arc-weighted  $f'$-degenerate. If $g'=0$, then $g=1_{\{x\}}$.  
  Let $X=\{x\}$. From the definition it follows that $G[X]$ is arc-weighted $(g-1)$-degenerate and $G-X$ is arc-weighted $(f-g)$-degenerate. 

   Otherwise, assume $g' \neq 0$. By the induction hypothesis,  there is a non-empty  subset $X$ of $V(G)$, such that $G[X]$ is arc-weighted   $(g'-1)$-degenerate, and $G-X$ is arc-weighted   $(f'-g')$-degenerate. Then $G[X]$ is arc-weighted  $(g-1)$-degenerate, and $G-X$ is arc-weighted $(f-g)$-degenerate.

Assume $f(x) = d_{(D,w)}^+(x)$ for each vertex $x$. 
By definition, $(D,w)$ has a dominating arc $e=(x,y)$.  As $(D,w)$ is a  minimal arc-weighted   acyclic orientation,  $w(e) = d_{(D,w)}^+(y)+1 =f(y)+1$.
Then $G'=G-e$ is  arc-weighted   $f'$-degenerate, where $f'=f-(f(y)+1)1_{\{x\}}$.

    Let $g' \in \mathbb{N}^{G'}$ be defined as $g'(v)=g(v)$, except that $$g'(x)=\min\{f'(x)+1, \max\{0,g(x)-g(y)\}\}.$$ 
   If $g'=0$, then  $g(x) \le g(y)$ and $g(v)=0$ for $v \ne x$. But $g(x) \le g(y) = 0$ implies that $g=0$, a contradiction.
   Thus $g' \ne 0$. By the induction hypothesis,  there is a non-empty   subset $X$ of $V(G')$, so that 
$G'[X]$ is arc-weighted  $(g'-1)$-degenerate, and $G'-X$ is arc-weighted   $(f'-g')$-degenerate. 

\bigskip 
\noindent
{\bf Case 1}  $x \in X$.

As $G-X=G'-X$ is arc-weighted   $(f'-g')$-degenerate, and $(f-g)(v)= (f'-g')(v)$ for $v \in V(G-X)$,   we conclude that  $G-X $ is arc-weighted  $(f-g)$-degenerate. 
It remains to show that $G[X]$ is arc-weighted $(g-1)$-degenerate.

If $y \notin X$, then
$G[X]=G'[X]$ is arc-weighted  $(g'-1)$-degenerate, and hence arc-weighted $(g-1)$-degenerate, since $g'(v)\le g(v)$ for every $v\in X$.

Assume that $y \in X$.
 As $G'[X]$ is arc-weighted   $(g'-1)$-degenerate, and $x \in X$, we have   $g'(x) \ge 1$.  Hence 
$g'(x) = \min\{f'(x)+1, g(x)-g(y)\} \le g(x)-g(y)$.  Let $(D',w')$ be  an arc-weighted   acyclic orientation  of $G'[X]$ with $d_{(D',w')}^+  \le g'-1$. 
 Let $(D'',w'')$ be obtained from $(D',w')$ by adding the arc $e=(x,y)$ of weight $w''(e)= g(y)=g'(y)$. Then $$w''(e)   = g'(y)   \ge d_{(D',w')}^+(y)+1 = d_{(D'',w'')}^+(y)+1,$$
and $$d_{(D'',w'')}^+(x) = d_{(D',w')}^+(x)+w''(e) < g'(x)+g(y) \le  g(x).$$
Hence, $e$ is a dominating arc in $(D'',w'')$ and $(D''-e,w'') = (D',w')$. Therefore $(D'',w'')$ is an arc-weighted acyclic orientation of $G[X]$ with 
$ d_{(D'',w'')}^+  \le  g-1$. Hence $G[X]$ is arc-weighted   $(g-1)$-degenerate.

\bigskip 
\noindent
{\bf Case 2}   $x \notin X$. 

Then  $G[X] = G'[X]$ is arc-weighted  $(g'-1)$-degenerate,  and $g'(v)=g(v)$ for all $v \in X$. Hence $G[X]$ is arc-weighted   $(g-1)$-degenerate. 
It remains to show that $G-X$ is arc-weighted $(f-g)$-degenerate.

As $G'-X$ is arc-weighted $(f'-g')$-degenerate,  we know that $(f'-g')(x) \ge 0$. So $g'(x) \ne f'(x)+1$ and hence $$ g'(x)= \max\{0, g(x)-g(y)\} \ge g(x)-g(y).$$ 
This implies that $(f'-g')(x) \le (f(x)-(f(y)+1)) - (g(x)-g(y)) \le (f-g)(x)$   (note that $f(y)+1 \ge g(y)$). 
If   $y \in X$,   then   $G-X = G'-X$ is arc-weighted  $(f-g)$-degenerate. 
  Therefore, we assume that $y \notin X$. 

  Let $(D',w')$ be an arc-weighted  acyclic orientation of $G'-X$ with $d_{(D',w')}^+   \le (f'-g')$.  
  Let $(D'',w'')$ be obtained from $(D',w')$ by adding the arc $e=(x,y)$ of weight $w''(e)=f'(y)-g'(y)+1=f(y)-g(y)+1$. Thus, $(D'',w'')$ is an arc-weighted   acyclic orientation of $G-X$. As $d_{(D'',w'')}^+(x)   \le (f'-g')(x) + w''(e) \le (f-g)(x)$, we have  $d_{(D'',w'')}^+  \le f-g$, and hence $G-X$ is arc-weighted   $(f-g)$-degenerate.
 
   This completes the proof of Lemma \ref{lem-1}.
		\end{proof}

\begin{theorem}
\label{thm-DP}
If $G$ is arc-weighted  $f$-degenerate, then $G$ is   $(f+1)$-DP-paintable.  
\end{theorem}
\begin{proof}
  We prove this by induction on the number of vertices. If $V(G)=\emptyset$, then this follows from the definition.
    Assume $V(G) \ne \emptyset$. Assume $g \in \mathbb{N}^G$ and $0 \ne g \le f+1$ and $(L,M)$ is a $g$-DP-cover of $G$. We shall show that there is a non-empty subset $X$ of $V(G)$ such that $G[X]$ is $g$-DP-colourable, and $G-X$ is $(f+1-g)$-DP-paintable.  
   By Lemma \ref{lem-1},  there is a non-empty subset $X$ of $V(G)$ such that $G[X]$ is arc-weighted  $(g-1)$-degenerate, and $G-X$ is arc-weighted  $(f-g)$-degenerate. According 
to Lemma \ref{lem-dpf}, 
    $G[X]$ is $g$-DP-colourable, and according to the induction hypothesis, $G-X$ is $(f+1-g)$-DP-paintable. By Lemma \ref{lem-0}, $G$ is $(f+1)$-DP-paintable.
\end{proof}

\section{Arc-weighted acyclic orientations of planar graphs}
\label{degree-truncated}

Thomassen \cite{Thomassen} proved that every planar graph is $5$-choosable. 
This result was strengthened in \cite{DP} to show that every planar graph is $5$-DP-colourable.
Hefetz \cite{Hefetz} asked the question whether every planar graph is $5$-AT, and the question remained open for a decade, and 
an affirmative answer was given in \cite{Zhu5AT}, and two alternate proofs of this result were given in \cite{GZ} and \cite{kozik}. The proofs in \cite{DP,Zhu5AT,GZ,kozik} are all parallel to Thomassen's proof. This is not a coincidence. The following result shows that Thomassen's proof can be easily adapted to show 
 that every planar graph is arc-weighted $4$-degenerate and hence $5$-DP-paintable, as well as $5$-AT.

\begin{theorem}
    \label{thm-planar}
    Assume $G$ is a plane graph with boundary $B(G)=[v_1v_2\ldots v_n]$.
    Let \[
    f_{G, v_1v_2} (v) = \begin{cases} 0, &\text{ if $v \in \{v_1,v_2\}$}, \cr
    2, &\text{ if $v \in B(G) - \{v_1,v_2\}$,} \cr 
    4, &\text{ if $v$ is an interior vertex}.   
    \end{cases}
    \]
   Then $G-v_1v_2$ is arc-weighted   $f_{G,v_1v_2}$-degenerate.
\end{theorem}
\begin{proof}
    The proof is by induction on $|V(G)|$. The case $|V(G)| \le 3$ is trivial. Assume $|V(G)|\ge 4$, and the theorem holds for graphs with fewer vertices. If $B(G)$ has a chord $e=v_iv_j$ that separates $G$ into two parts $G_1,G_2$, with $v_1, v_2 \in V(G_1)$, then we apply the induction hypothesis to obtain an arc-weighted   acyclic orientation $(D_1, w_1)$ of $G_1-v_1v_2$ with $d_{(D_1,w_1)}^+(v) \le f_{G_1, v_1v_2}(v)$ and an arc-weighted    acyclic orientation $(D_2, w_2)$ of $G_2-v_iv_j$ with $d_{(D_2,w_2)}^+(v) \le f_{G_2, v_iv_j}(v)$.
    
        Now we show that  
    $(D,w)= (D_1 \cup D_2, w_1\cup w_2)$ is an arc-weighted  acyclic orientation  of $G-v_1v_2$ with $d_{(D,w)}^+(v) \le f_{G, v_1v_2}(v)$. We extend $f_{G_1,v_1v_2}$ to $V(G)$ by letting $f_{G_1,v_1v_2}(v)=0$ for $v \in V(G)-V(G_1)$, and extend $f_{G_2,v_iv_j}$ to $V(G)$ similarly. Then 
   $f_{G_1,v_1v_2}+f_{G_2, v_iv_j} \le  f_{G, v_1v_2}$. So $d_{(D,w)}^+(v) \le f_{G, v_1v_2}(v)$ for each vertex $v$ of $G$. 

   For any non-empty sub-digraph $D'$ of $D$, if $E(D') \cap E(D_1) \ne \emptyset$, then since $(D_1, w_1)$ is acyclic, there is an  arc $e=(x,y) \in E(D') \cap E(D_1)$ that is dominating in $(D' \cap D_1, w)$. As the arcs in $D_2$ do not increase the out-degree of $y$,  $e$ is also a dominating arc for $D'$. If $E(D') \cap E(D_1) = \emptyset$, then  there is an arc $e=(x,y) \in E(D') \cap E(D_2)$ that is a dominating arc in $( D' \cap D_2, w)$, which is also a dominating arc in $(D',w)$.

    Assume $B(G)$ has no chord. Let $G'=G-v_n$. Let $u_1,u_2,\ldots, u_k$ be the interior neighbours of $v_n$. By the induction hypothesis, there is an arc-weighted    acyclic orientation $(D', w')$ of $G'-v_1v_2$ with $d_{(D',w')}^+(v) \le f_{G', v_1v_2}(v)$ for all $v \in V(G')$. We extend $D'$ to an arc-weighted orientation $(D,w)$ of $G$, by adding arcs $(v_n, v_1), (v_n,v_{n-1})$ of weight $1$, and arcs $(u_i,v_n)$ for $i=1,2,\ldots, k$  of weight $2$. 

 It follows from the definition that $d_{(D,w)}^+(v) \le f_{G, v_1v_2}(v)$ for all $v \in V(G)$.
    Now we show that $(D,w)$ is    an arc-weighted   acyclic  orientation of $G-v_1v_2$. Let $D'$ be a non-empty sub-digraph of $D$. We  show that $(D',w)$ has a dominating arc.
    \begin{enumerate}
        \item If $(v_n, v_1) \in E(D')$, then $(v_n, v_1)$ is a dominating arc in $(D',w)$. 
        \item Otherwise,  if 
    $(u_i, v_n) \in E(D')$ for some $i$, then $(u_i, v_n)$ is a dominating arc in $(D',w)$.
    \item Otherwise, if $E(D') - \{(v_n, v_{n-1})\} \ne \emptyset$, then $(D'- \{(v_n, v_{n-1})\},w)$ has a dominating arc by the induction hypothesis, which is also a dominating arc in $(D',w)$.
    \item Otherwise, $(v_n, v_{n-1})$ is a dominating arc of $D'$. \qedhere
    \end{enumerate}
\end{proof}

\section{Degree-truncated arc-weighted strict degeneracy}

A graph $G$ is called {\em degree-choosable} (respectively, {\em degree-DP-colourable}) if $G$ is $d_G$-choosable (respectively, $d_G$-DP-colourable). 

The class of degree-choosable graphs was characterized in \cite{Borodin, ERT}: A connected graph $G$ is not degree-choosable if and only if it is a Gallai-tree, i.e., each block of $G$ is either a clique or an odd cycle (a block of $G$ is a maximal 2-connected subgraph of $G$). The class of degree-DP-colourable graphs was characterized in \cite{BKP,KO19}: A connected graph $G$ is not degree-DP-colourable if and only if it is a GDP tree, i.e., each block of $G$ is either a clique or a  cycle.

 Lemma \ref{lem-degree} below characterizes arc-weighted $(d_G-1)$-degenerate graphs.
   
\begin{lem}\label{lem-degree}
     A connected graph $G$ is not arc-weighted $(d_G-1)$-degenerate if and only if $G$ is a GDP-tree.
\end{lem}
\begin{proof}
If $G$ is a GDP-tree, then $G$ is not degree-DP-colourable, and hence it is not arc-weighted $(d_G-1)$-degenerate.

Assume that $G$ is a connected graph and is not a GDP-tree. Then $G$ has a block $B$ that is   neither a cycle nor a complete graph. We claim that $B$ has an induced theta subgraph $H$, i.e., three internally disjoint paths connecting two vertices. Indeed, if $B$ is chordal, since $B$ is not complete, $B$ contains at least two maximal cliques. As $B$ is 2-connected, the two maximal cliques have at least two common vertices. Choosing two vertices in their intersection and one vertex from each clique outside the intersection gives $H$. 
If $B$ is not chordal, then $B$ has an induced cycle $C$ of length at least $4$. If there is a vertex $v \in V(B)-V(C)$ adjacent to at least three vertices of $C$, then $v$ together with a segment of $C$ induces a theta graph. Otherwise, there is an induced path $P$ connecting two vertices of $C$, whose interior vertices are disjoint from $C$ and not adjacent to the vertices in $C$. Then $P \cup C$ is again an induced theta graph. 

 Let $H$ be 
 an induced theta graph, and let   $e=uv$ be an edge of  $H$ with $d_H(u)=3$ and $  d_H(v)=2$. Then $H$ has an acyclic orientation $D'$ 
    in which $v$ is the unique source and $u$ is the unique sink. Reverse the direction of the edge $uv$ as an arc $(u,v)$ of weight $2$, we obtain an arc-weighted digraph $(D,w)$  (other arcs of $D$ are of weight $1$). For any non-empty sub-digraph $D''$ of $D$, either $(u,v) \in E(D'')$ and hence is a dominating arc in $(D'', w)$ or $D''$ has an arc $(x,y)$ with $y$ being a sink and hence is a dominating arc. So $(D,w)$  
     is acyclic and $d_{(D,w)}^+(x) \le d_H(x)-1$ for each vertex $x$.
    Let $D''$ be an acyclic orientation of $G-E(H)$ such that the vertices of $H$ are the only sources. Such an orientation exists, for example, by ordering the vertices so that every vertex in $V(G)-V(H)$ has a neighbor preceding it and all vertices of $H$ come first and then orienting every edge from the earlier endpoint to the later endpoint.
    Then $(D \cup D'', w)$ (with all other arcs of weight 1) is acyclic and each vertex $x$ has 
    $d_{(D \cup D'', w)}^+(x) \le d_G(x)-1$. Thus, $G$ is arc-weighted $(d_G-1)$-degenerate.
\end{proof}

\begin{definition}
    \label{def-degree-truncated}
    Let $G$ be a graph, and $k$ be a positive integer. Let 
$$\phi_k(v)=\min\{k, d_G(v)\}.$$
If  $G$ is  $\phi_k$-choosable, then we   say  $G$ is  degree-truncated $k$-choosable. 
    The  degree-truncated choice number of $G$ is 
$$\operatorname{ch}^{\text{\st{d}}}(G)= \min\{k:\text{ $G$ is degree-truncated  $k$-choosable}\},$$ 
and for 
a family  $ \mathcal{G}$ of graphs, $$\operatorname{ch}^{\text{\st{d}}}(\mathcal{G}) = \max \{\operatorname{ch}^{\text{\st{d}}}(G): G \in \mathcal{G}\}.$$
\end{definition}

 If $G$ is not degree-choosable, then $G$ is not degree-truncated $k$-choosable for any $k$. Thus, the degree-truncated choice number of Gallai-trees is not defined. For connected graphs $G$ other than Gallai-trees, $\operatorname{ch}^{\text{\st{d}}}(G)$ is well-defined and is at most the maximum degree of $G$.

Problems concerning the degree-truncated choice number of planar graphs was first considered by Bruce Richter (see \cite{Hutchinson}). For any positive integer $k$, the complete bipartite graph $K_{2, k^2}$ is planar, not a Gallai-tree, and not degree-truncated $k$-choosable. Thus there is no constant upper bound on the degree-truncated choice number of planar graphs that are not Gallai-trees.  

Let $\mathcal{P}$ be the family of 3-connected non-complete planar graphs. Richter \cite{Hutchinson} asked whether $\operatorname{ch}^{\text{\st{d}}}(\mathcal{P}) \le 6$. 

Richter's question has attracted some attention \cite{Hutchinson,CPTV} and remained open for a decade. It also remained open whether there is a constant upper bound for $\operatorname{ch}^{\text{\st{d}}}(\mathcal{P})$. Recently,  a  negative answer was given in \cite{ZZZ25}, where a graph $G \in \mathcal{P}$ with $\operatorname{ch}^{\text{\st{d}}}(G) = 8$ was constructed. In \cite{ZZZ25} it was also proved that $\operatorname{ch}^{\text{\st{d}}}(\mathcal{P}) \le 16$.
These results are improved in \cite{JXXZ} and \cite{XZZZ}. In \cite{XZZZ}, a graph $G \in \mathcal{P}$ with $\operatorname{ch}^{\text{\st{d}}}(G) = 10$ was given, and it was also proved in   \cite{XZZZ} that  
$\operatorname{ch}^{\text{\st{d}}}(\mathcal{P}) \le 11$. 
So $\operatorname{ch}^{\text{\st{d}}}(\mathcal{P})=10$ or $11$, and it was conjectured in \cite{XZZZ}
that $\operatorname{ch}^{\text{\st{d}}}(\mathcal{P})=10$.

We generalize this concept to the 
 {degree-truncated arc-weighted strict degeneracy} of graphs. 
 
 \begin{definition}
     \label{def-degree-truncated-st-deg}
     Let $G$ be a graph, and $k$ be a positive integer. Let $\phi_k(v) = \min\{k, d_G(v)\}$. If   
     $G$ is arc-weighted  $(\phi_k-1)$-degenerate, then we say 
       $G$ is {\em degree-truncated arc-weighted strict $k$-degenerate}.  
    The degree-truncated arc-weighted strict degeneracy of $G$ is defined as  $$\operatorname{sd}_{w}^{\text{\st{d}}}(G)= \min \{k: \text{$G$ is arc-weighted  $(\phi_k-1)$-degenerate}\}.$$ 
    For a family $\mathcal{G}$ of graphs, 
    $$\operatorname{sd}_w^{\text{\st{d}}}(\mathcal{G})=\max \{ \operatorname{sd}_w^{\text{\st{d}}}(G): G \in \mathcal{G}\}.$$
\end{definition}

Given a graph $G$, the degree-truncated DP-chromatic number $\chi_{\mathrm{DP}}^{\text{\st{d}}}(G)$ of $G$, the degree-truncated DP-paint number $\chi_{\mathrm{DPP}}^{\text{\st{d}}}(G)$ of $G$, the degree-truncated paint number $\chi_{P}^{\text{\st{d}}}(G)$ of $G$ and the degree-truncated AT number $\operatorname{AT}^{\text{\st{d}}}(G)$ of $G$ are defined similarly. For a family $\mathcal{G}$ of graphs,  
$\chi_{\mathrm{DP}}^{\text{\st{d}}}(\mathcal{G})$, 
$\chi_{\mathrm{DPP}}^{\text{\st{d}}}(\mathcal{G})$,
$\chi_{P}^{\text{\st{d}}}(\mathcal{G})$, and 
$\operatorname{AT}^{\text{\st{d}}}(\mathcal{G})$ are the maximum of the corresponding parameters for the graphs in $\mathcal{G}$. It follows from Theorems \ref{thm-AT} and \ref{thm-DP} that for any family $\mathcal{G}$ of graphs, 
$$\chi_{\mathrm{DP}}^{\text{\st{d}}}(\mathcal{G}) \le \chi_{\mathrm{DPP}}^{\text{\st{d}}}(\mathcal{G}) \le \operatorname{sd}_w^{\text{\st{d}}}(\mathcal{G}),  \text{ and }   \chi_{P}^{\text{\st{d}}}(\mathcal{G}) \le \operatorname{AT}^{\text{\st{d}}}(\mathcal{G}) \le \operatorname{sd}_w^{\text{\st{d}}}(\mathcal{G}).$$

 In \cite{ZZZ25} it was proved that $\chi_{\mathrm{DP}}^{\text{\st{d}}}(\mathcal{P}) \le 16.$
 Theorem \ref{thm-stwdp} below implies that $$\chi_{\mathrm{DP}}^{\text{\st{d}}}(\mathcal{P}), \, \chi_{\mathrm{DPP}}^{\text{\st{d}}}(\mathcal{P}),   \, \chi_{P}^{\text{\st{d}}}(\mathcal{P}),\,  \operatorname{AT}^{\text{\st{d}}}(\mathcal{P}) \le 15.$$

\begin{theorem}
    \label{thm-stwdp}
    For the family  $  \mathcal{P}$ of 3-connected non-complete planar graphs, we have
     $\operatorname{sd}_w^{\text{\st{d}}}(\mathcal{P})  \le 15.$ 
\end{theorem}  
\begin{proof}
Assume that the theorem is not true and $G$ is a counterexample with minimal number of vertices and subject to this, with maximum  number of edges.

Let $$V_1=\{v \in V(G): d_G(v) < 15\} \text{ and  } V_2 = \{v \in V(G): d_G(v) \ge 15\}.$$
Let $\phi(v)= \min \{15, d_G(v)\}.$

If $V_2=\emptyset$, then $\phi=d_G$. Since a $3$-connected
non-complete graph is not a GDP-tree, Lemma~\ref{lem-degree} implies
that $G$ is arc-weighted $(\phi-1)$-degenerate. Hence, we may assume
that $V_2\ne\emptyset$.

Let $\mathcal{Q}$ be the set of connected components of $G[V_1]$.
For each $Q\in\mathcal{Q}$, let $\theta_Q$ be the face of   $G[V_2]$ that contains $Q$, and let $V(\theta_Q)$ denote the set of vertices on the boundary of $\theta_Q$. Note that $G[V_2]$ need not be connected, and hence $G[V(\theta_Q)]$ need not be connected. 

Note that any two vertices $u,v$ of $V_2$ on the boundary of a same face of $G$ are adjacent, for otherwise adding the edge $uv$ is a counterexample to Theorem \ref{thm-stwdp} with more edges than $G$, in contradiction to the choice of $G$.

 As $G$ is $3$-connected, this implies that 
each face $\theta$ of $G[V_2]$ contains at most one connected component of
$G[V_1]$. So the faces $\theta_Q$, $Q\in\mathcal{Q}$, are pairwise distinct. Also  
every vertex $v$ of $V(\theta_Q)$ has a neighbor in $Q$, for otherwise, by adding an edge joining the two neighbors of $v$ in $V(\theta_Q)$ results in a counterexample with more edges.

Given a non-empty plane graph $X$, let $\Theta(X)$ be the bipartite graph with partite sets $V(X)$ and $F(X)$ (the set of faces of $X$) in which $v\theta$ is an edge if and only if $v \in V(\theta)$, i.e., $v$ is on the boundary of $\theta$.

An acyclic orientation of a graph $H$  with exactly one source and one sink is called a {\em  bipolar orientation} of $H$. It is well-known \cite{LEC1967, FOR1995} that for a 2-connected plane graph $H$ and two adjacent vertices $u,v$ on the boundary of the infinite face of $H$, there is  a bipolar orientation $D$ such that 
\begin{enumerate}
    \item $u$ is the unique source, $v$ is the unique sink, and
    \item for   each face $\theta$ of $H$, the induced sub-digraph $H[V(\theta)]$ also has a unique source $s_{\theta}$ and a unique sink $t_{\theta}$, and 
    \item for each vertex $w \in V(H)-\{u,v\}$, 
    there are  exactly two faces $\theta$ incident to $v$ such that $v \notin \{s_\theta, t_\theta\}$. 
\end{enumerate}   

\begin{lem}
    \label{lem:assignment}
    Assume $H$ is a non-empty plane graph and $uv$ is an edge on the boundary of the infinite face of $H$. There exists,    for each   face $\theta $, a subset $P(\theta)$ of $V(\theta)$ such that the following hold:
    \begin{enumerate}
    \item For each finite face $\theta \in F(H)$, 
$|P(\theta)| =|V(\theta)|-2 $, and for the infinite face $\theta^*$, $P(\theta) =V(\theta^*)$.
\item Let
$V'$ be the set of vertices of $H$ incident to only the infinite face.
Then each vertex in $V'\cup \{u,v\}$ is contained only in $P(\theta^*)$, and each other
vertex  of $H$ is contained in at most two of the sets $\{P(\theta): \theta \in F(H)\}$.
\item If  $v_1,v_2 \in V(\theta)-P(\theta)$, then there is a cycle $C$ in $H$ containing $v_1$ and $v_2$ such that $\theta$ is contained in the interior of $C$.
\end{enumerate}
\end{lem} 
\begin{proof}
 If $H$ is $2$-connected, then let $D$ be a bipolar orientation of $H$ with $u$ be the unique source and $v$ be the unique sink. 
For each interior face $\theta$ of $H$, let $P(\theta) =V(\theta)-\{s_{\theta}, t_{\theta}\}$ and for the infinite face $\theta^*$ of $H$, let $P(\theta) =V(\theta)$. Then each
of $u,v$ is contained only in $P(\theta^*)$, and each 
vertex $w \ne u,v $ of $H$ is contained in exactly two $P(\theta)$ for which $w \in V(\theta)$ and $ w \notin \{s_\theta, t_\theta\}$.

If $H$ is connected but not 2-connected, then let $B_1, B_2, \ldots, B_k$ be the blocks of $H$. Assume $uv$ is a boundary edge of $B_1$. For each $i \ge 2$, let $u_iv_i$ be a boundary edge of $B_i$ such that $u_i$ is the cut-vertex of $H$ contained in $B_i$ along the shortest path from $B_i$ to $B_1$. For each $i$, for each face $\theta \in F(B_i)$, let $P_{B_i}(\theta)$ be defined as above. 

Then for each face $\theta$ of $H$, 
there is a subset $I_\theta$ of $\{1,2,\ldots, k\}$ such that $\theta$ is the intersection $\cap_{i \in I_\theta} \theta_i$, where $\theta_i$ is a face of $B_i$. Hence $V(\theta) = \cup_{i \in I_\theta}V(\theta_i)$. Moreover, if $\theta$ is the infinite face of $H$, then  $\theta_i$ is the infinite  face of $B_i$ for each $i \in I_\theta$. Otherwise, there is exactly one index $j \in I_\theta$ such that $\theta_j$ is a finite face of $B_j$.   Let $P(\theta)=\cup_{i \in I_\theta}P_{B_i}(\theta_i)$. It is easy to verify that $P$ satisfies (1), (2) and (3).

If $H$ is not connected, then let $J_1, \ldots, J_k$ be the connected components of $H$. Similarly, for each face $\theta$ of $H$, there is a subset $I_\theta$ of $\{1,2,\ldots, k\}$ such that $\theta=\cap_{i \in I_\theta} \theta_i$ and $V(\theta) = \cup_{i \in I_\theta}V(\theta_i)$, where $\theta_i$ is a face of $J_i$. If $\theta$ is the infinite face of $H$, then  $\theta_i$ is the infinite  face of $J_i$ for each $i \in I_\theta$. Otherwise, there is exactly one index $j \in I_\theta$ such that $\theta_j$ is a finite face of $J_j$.   Let $P(\theta)=\cup_{i \in I_\theta}P_{J_i}(\theta_i)$. It is easy to verify that $P$ satisfies (1), (2) and (3).
\end{proof}

We apply Lemma~\ref{lem:assignment} to $G[V_2]$.
For $Q\in\mathcal{Q}$, vertices in $P(\theta_Q)$ are {\em  protectors} of $Q$, and vertices in $V(\theta_Q)-P(\theta_Q)$ are {\em non-protectors} of $Q$. Thus each
component $Q\in\mathcal{Q}$ has at most two non-protectors, and 
each vertex $v\in V_2$ is a protector of at most two
components of $\mathcal{Q}$.

\begin{lem}
    \label{lem:protector}
    For each component $Q\in\mathcal Q$, there is a leaf-block $B$ of $Q$ such that $B$ contains a non-root vertex adjacent to a protector of $Q$.
\end{lem}

\begin{proof}
Assume to the contrary that no protector of $Q$ is adjacent to a non-root vertex of a leaf-block of $Q$.

First, note that $Q$ has at least two leaf-blocks. Otherwise, $Q$ itself is a leaf-block, and every vertex of $Q$ is a non-root vertex. Recall that every vertex of $V(\theta_Q)$ is adjacent to some vertex of $B$. Thus at least one protector in $V(\theta_Q)$ is adjacent to $Q$, which contradicts the assumption.

Let $B_1$ and $B_2$ be two leaf-blocks of $Q$. By the assumption, we have $
N_{G[V_2]}(U(B_i))\subseteq \{v_1,v_2\}$ for $i=1,2$,
where $U(B_i)$ is the set  $V(B_i)$ with its unique cut vertex removed, and $v_1,v_2$ are the two non-protectors of $Q$ (if they exist).

If $N_{G[V_2]}(U(B_i))=\{v_1,v_2\}$ for $i=1,2$, then there is a cycle $C$ contained in
$\{v_1,v_2\}\cup U(B_1)\cup U(B_2)$
such that all vertices of $Q-C$ lie in the interior of $C$, while
all vertices of $V(\theta_Q)-\{v_1,v_2\}$ lie in the exterior of $C$. Hence no vertex of $V(\theta_Q)-\{v_1,v_2\}$ can be adjacent to $U(B_1)\cup U(B_2)$, contradicting the fact that every vertex of $V(\theta_Q)$ is adjacent to some vertex of $Q$.

Therefore, without loss of generality, assume that $v_2\notin N_G(U(B_1))$.
Let $x$ be the root of $B_1$. Since no protector of $Q$ is adjacent to a vertex of $U(B_1)$ and $v_1,v_2$ are the only non-protectors of $Q$, we have
$N_G(U(B_1))\setminus U(B_1)\subseteq\{x,v_1\}$. Since $Q$ has another leaf block, the set $\{x,v_1\}$ separates
$U(B_1)$ from the rest of $G$, contradicting the $3$-connectivity
of $G$.
\end{proof}

For each $Q \in \mathcal{Q}$, choose a non-root vertex $v_Q$ of a leaf block of $Q$ such that 
\begin{enumerate}
    \item $v_Q$ is adjacent to a protector of $Q$.
    \item Subject to (1), the number of non-protectors adjacent to $v_Q$ is minimum.
\end{enumerate}

For each $Q\in\mathcal Q$, we now describe how to orient the edges
inside $Q$ and the edges between $Q$ and $V_2$. 

Assume that  $Q$ is not a GDP-tree. Take an arc-weighted acyclic orientation $(D_Q,w_Q)$ of $Q$ given by Lemma~\ref{lem-degree}, such that
$d^+_{(D_Q,w_Q)}(v)\le d_Q(v)-1$ for every vertex $v\in V(Q)$.
Orient all edges between $Q$ and $V_2$ from $Q$ to $V_2$ with
weight $1$.

Assume that  $Q$ is a GDP-tree. Note that $d_Q(v_Q) \le 3$ and if $d_Q(v_Q)=3$, then the leaf-block of $Q$ containing $v_Q$ is a copy of $K_4$. We consider two cases. Arcs defined below have weight 1, unless it is specified otherwise.

\medskip
\noindent
{\bf Case 1} $d_Q(v_Q) \le 2$, or $d_Q(v_Q)=3$ and $v_Q$ is adjacent to at most one non-protector of $Q$.

In this case, we orient edges in $Q$ so that it is acyclic with $v_Q$ as the unique source, orient all edges between $Q$ and $V_2$ as arcs from $Q$ to $V_2$, except that if $v$ is a protector of $Q$ adjacent to $v_Q$, then the edge $vv_Q$ is oriented as an arc $(v, v_Q)$ of weight $5$.

\medskip
\noindent
{\bf Case 2}  $d_Q(v_Q)=3$ and $v_Q$ is adjacent to two non-protectors of $Q$.

In this case,  the leaf-block $B$ of $Q$ containing $v_Q$ is a copy of $K_4$, and by Lemma \ref{lem:assignment}, there is a cycle $C$ contained in $V(\theta_Q)$ containing the two non-protectors of $Q$, and $Q$ is contained in the interior of $C$. Then $B$ has one vertex $u_Q$ that is  contained in the interior of the triangle $T$ induced by the other three vertices of $B$. Note that since $v_Q$ is adjacent to the two non-protectors in $C$, $v_Q \ne u_Q$.

If $u_Q$ is the root of $B$ (and hence a cut-vertex of $Q$), then there is another leaf-block $B'$ 
contained in the 
interior of $T$. Then a non-root vertex $v'$ of $B'$ is adjacent to a vertex of $V_2$. As the two non-protectors 
of $Q$ lie in the exterior of $T$, $v'$ is not adjacent to any non-protector of $Q$. Hence $v'$ is adjacent to a protector of $Q$. This contradicts the choice of $v_Q$. Thus $u_Q$ is a non-root vertex of $B$. If $u_Q$ is adjacent to some vertex of $V_2$, then $u_Q$ is adjacent to a protector of $Q$, but not adjacent to any non-protector of $Q$. Again 
this contradicts the choice of $v_Q$. Thus $u_Q$ is not adjacent to any vertex of $V_2$.
We orient all edges between $Q$ and $V_2$ as arcs from $Q$ to $V_2$, except that if $v$ is a protector of $Q$ adjacent to $v_Q$, then the edge $vv_Q$ is oriented as an arc $(v, v_Q)$ of weight $5$. 
Orient edges in $Q$ so that it is acyclic with $v_Q$ as the unique source, and with $u_Q$ be a sink.

\medskip

After this process has been completed for all $Q \in \mathcal{Q}$, we obtain an arc-weighted orientation $(D',w')$ of $G-E(G[V_2])$. Let $(D'',w'')$ be an arc-weighted acyclic orientation of $G[V_2]$ with $d_{(D'',w'')}^+(v) \le 4$ for each vertex $v \in V_2$. By Theorem \ref{thm-planar}, such an orientation exists. Let $(D,w) =(D',w') \cup (D'',w'')$. 
Now we show that $(D,w)$ is an arc-weighted acyclic orientation and $d_{(D,w)}^+(v) \le \phi(v)-1$ for each vertex $v$. 

First, we check the out-degree of the vertices of $D$. If $Q$ is not a GDP-tree, then for every $v\in V(Q)$,
$d^+_{(D,w)}(v)
 \le d^+_{(D_Q,w_Q)}(v)+d_G(v)-d_Q(v)
 \le d_G(v)-1$.
If $Q$ is a GDP-tree, then all outgoing arcs of a vertex $v\in V(Q)$ have weight $1$, and every such vertex has an incoming arc. Hence $d^+_{(D,w)}(v)\le d_G(v)-1$. For $v \in V_2$, an edge $v w$ with $w\in V_1$ is an out-going arc of $v$ if and only if $w=v_Q$ for some $Q \in \mathcal{Q}$ such that $v$ is a protector of $Q$. As $v$ is a protector of at most two components $Q \in \mathcal{Q}$, and each such arc $(v, v_Q)$ has weight $5$, we know that such arcs contribute at most $10$ to $d_{(D,w)}^+(v)$. As $d_{(D'',w'')}^+(v) \le 4$, we conclude that  $d_{(D,w)}^+(v) \le 14 = \phi(v)-1$.

We now show that $(D,w)$ is arc-weighted acyclic. Let $H$ be an arbitrary non-empty subdigraph of $D$. For each $Q\in\mathcal Q$, let $E_Q$ be the set of arcs of $D$ incident with a vertex of $Q$,
and let
$$\mathcal Q_H=
\{Q\in\mathcal Q:E(H)\cap E_Q\ne\emptyset\}.$$
We prove by induction on $|\mathcal Q_H|$ that $H$ contains a dominating arc.

If $\mathcal Q_H=\emptyset$, then $H$ is a non-empty subdigraph of $D''$, and hence contains a dominating arc. Now assume that $\mathcal Q_H\ne\emptyset$, and choose $Q\in\mathcal Q_H$.

First suppose that $Q$ is not a GDP-tree. If $H-E_Q$ is non-empty, then, by induction, $H-E_Q$ contains a dominating arc. Since all arcs between $Q$ and $V_2$ are directed from $Q$ to $V_2$, restoring the arcs in $E_Q$ does not increase the out-degree of the head of this arc. Hence it remains dominating in $H$. If $H-E_Q$ is empty,
then either $H$ contains an arc from $Q$ to $V_2$, which is dominating since its head is a sink in $H$, or $H$ is a non-empty subdigraph of $(D_Q,w_Q)$ and hence contains a dominating arc.

It remains to consider the case that $Q$ is a GDP-tree. If $Q$ is oriented as in Case~2 and $H$ contains an arc incident with $u_Q$, then this arc is dominating, since every such arc is directed into
$u_Q$ and $u_Q$ has no neighbour in $V_2$.

Otherwise, suppose that $H$ contains an arc $(v,v_Q)$, where $v$ is a protector of $Q$. In Case~1, we have $d^+_{(H,w)}(v_Q)\le 4$.
The same inequality holds in Case~2, since the arc between $v_Q$ and $u_Q$ is not present in $H$. Therefore, $(v,v_Q)$ is dominating, as $w((v,v_Q))=5>d^+_{(H,w)}(v_Q)$.

Finally, suppose that neither of the preceding cases applies. If $H-E_Q$ is non-empty, then it contains a dominating arc by induction. Since $H$ contains no arc from a protector of $Q$ to $v_Q$, every arc of $H\cap E_Q$ between $Q$ and $V_2$ is directed from $Q$ to $V_2$. Thus restoring $E_Q$ does not increase the out-degree of the head of the dominating arc, so it remains dominating in $H$.

If $H-E_Q$ is empty, then either $H$ contains an arc from $Q$ to $V_2$, whose head is a sink in $H$, or $H$ is contained in the acyclic orientation of $Q$. In the latter case, an arc incident with a sink of $H$ is dominating. Thus, in every case, $H$ contains a dominating arc.

This completes the proof of Theorem \ref{thm-stwdp}.
\end{proof}

A graph $H$ is a {\em minor} of a graph $G$ if a copy of $H$ can be obtained from $G$ by deleting vertices and edges and contracting edges. A family $\mathcal{G}$ of graphs is {\em minor-closed} if $G \in \mathcal{G}$ and $H$ is a minor of $G$ implies that $H \in \mathcal{G}$. If $\mathcal{G}$ is minor-closed and does not contain all the graphs, then $\mathcal{G}$ is   {\em proper minor-closed}.  It was proved in \cite{ZZZ25} that if  $\mathcal{G}$ is a proper minor-closed family of graphs, and $G \in \mathcal{G}$ is highly connected and non-complete, then $\operatorname{ch}^{\text{\st{d}}}(G) \le k $ for some constant $k$. 
Theorem \ref{thm-minor} strengthens this result to arc-weighted strict degeneracy.  The proof uses a similar idea as in the proof of Theorem \ref{thm-stwdp}, and is  parallel to the proof of Theorem 5 in \cite{ZZZ25}.

\begin{theorem}
    \label{thm-minor}
    Assume $\mathcal{G}$ is a proper minor-closed family of graphs. Let $s$ be the minimum integer such that $K_{s,t} \notin \mathcal{G}$ for some positive integer $t$. Then there is a constant $k=k(s)$ such that every $s$-connected graph in $\mathcal{G}$ which is not a GDP-tree is degree-truncated arc-weighted strict $k$-degenerate.
    In particular, for any surface $\Sigma$, there is a constant $k$ such that every 3-connected non-complete graph embedded in $\Sigma$ is degree-truncated arc-weighted strict $k$-degenerate.
\end{theorem}
 \begin{proof}

First we note that for any proper minor-closed family $\mathcal{G}$ of graphs, there exist positive integers $s,t$ such that the complete bipartite graph $K_{s,t} \notin \mathcal{G}$. This is so, because there is a complete graph $K_s \notin \mathcal{G}$ (for otherwise   $K_s \in \mathcal{G}$ for all $s$ and hence every graph is in $\mathcal{G}$ and $\mathcal{G}$ is not proper), and $K_s$ is a minor of $K_{s,s}$. So the integer $s$ as described in the theorem is well-defined. 

Let $t$ be the minimum integer such that $K_{s,t} \notin \mathcal{G}$. 
Let $$q = 4^{s+1}s! st (s+t-1)+1, \text{ and }  k=2^{s+2}tq.$$  Let $$f(v)=\min \{k-1, d_G(v)-1\}.$$
We prove that $G$ is arc-weighted  $f$-degenerate by using
the alternate definition of 
arc-weighted $f$-degeneracy given in Lemma \ref{lem-alternate}. For simplicity, we also use the vertex-deletion operation (recall that $\dd_v$ can be realized by a sequence of arc-deletion operations, as described in Section~1). Thus we shall apply the two operations: arc-deletion and vertex-deletion to $(G,f)$ to obtain a no-vertex graph.

If $s=1$, then graphs in $\mathcal{G}$ have maximum degree at most $t-1$.  For $G \in \mathcal{G}$ which is not a GDP-tree, $G$ is arc-weighted  $(d_{G}-1)$-degenerate, and hence degree-truncated arc-weighted $(k-1)$-degenerate. Assume $s \ge 2$.

Let $V_1 = \{v \in V(G): d(v) < k\}$ and $V_2 =V(G)-V_1$. Let $p=|V_2|$. Using the fact that $G$ is $K_{s,t}$-minor-free,  we can show that $G$ is $(2^{s+2}t-1)$-degenerate.

The following two lemmas, proved in~\cite{Thomason07}, are needed
in the proof.

\begin{lem}
	\label{lem-1b} Every non-empty graph $G$ with at least $2^{s+1}t|V(G)|$ edges has a $K_{s,t}$-minor.
\end{lem}

\begin{lem}
	\label{lem-2}
If $G$ is a bipartite graph with partite sets $A$ and $B$, and with at least $(s-1)|A|+4^{s+1}s!t|B|$ edges, then $G$ has a $K_{s,t}$-minor.
\end{lem}

The following is a consequence of  Lemma \ref{lem-2}.

\begin{cor}
	\label{cor-1}
	Assume $G$ is a bipartite graph with partite sets $A$ and $B$, and each vertex in $A$ has degree at least $s$. If $G$ has no $K_{s,t}$-minor, then $|E(G)| < 4^{s+1}s! st|B|$ and consequently, $B$ has a vertex $v$ of degree at most $4^{s+1}s! st$.
\end{cor}
\begin{proof}
	Assume $G$ has no $K_{s,t}$-minor. By Lemma \ref{lem-2}, $|E(G)| < (s-1)|A|+4^{s+1}s!t|B|$. 
	So $|E(G)| - (s-1)|A| < 4^{s+1}s!t|B|$. 
	
	As each vertex in $A$ has degree at least $s$, we have $|E(G)| \ge s|A|$, which implies that $|E(G)| - (s-1)|A| \ge |E(G)|- \frac{s-1}{s}|E(G)| = |E(G)|/s$. Therefore 
	$|E(G)|  \le s (|E(G)| - (s-1)|A|) < 4^{s+1}s! st|B|$.
\end{proof}

It follows from Lemma \ref{lem-1b} that each subgraph of $G$ has minimum degree at most $2^{s+2}t-1$, that is, $G$ is $(2^{s+2}t-1)$-degenerate. Hence the vertices of $V_2$ can be ordered as 
  $w_1,w_2, \ldots, w_p$  so that each $w_i$ has at most $2^{s+2}t-1$ neighbors $w_j$ with $j < i$.

We process the vertices of $V_2$ in the order
$w_1,w_2,\ldots,w_p$. When processing $w_i$, for each neighbor $w_j$ of $w_i$ with $j>i$, we apply the arc-deletion operation
$$
\ed_{(w_j,w_i),\,f(w_i)-q+1},
$$
where $f(w_i)$ denotes the current value of $f$ at the time the operation is applied. Proceeding in this way for
$i=1,2,\ldots,p$, all edges of $G[V_2]$ are deleted.
The choice of $k$ and the fact that $G$ is $(2^{s+2}t-1)$-degenerate guarantee that
all these operations are legal. Let $(G',f')$ be the resulting
pair. Then $V_2$ is an independent set in $G'$ and
$f'(w_i)\ge q-1$ for every $w_i\in V_2$.

We contract each connected component $Q$ of  $G'[V_1]$ into a vertex $v_Q$. Denote the resulting bipartite graph by  $H$. It follows from Corollary \ref{cor-1} that there is an ordering
  $u_1,u_2, \ldots, u_p$ of vertices of $V_2$, so that $u_i$ has degree at most  $ 4^{s+1}s!st$ in $H-\{u_p, u_{p-1}, \ldots, u_{i+1}\}$.

We process the vertices of $V_2$ one by one in the order $u_1,u_2, \ldots, u_p$.  Denote by $X$ the set of deleted vertices. Initially $X = \emptyset$. In the process, deleted vertices are added to $X$. So $X$ is a dynamic set, whose size increases in the process.

To process $u_i$, we do the following: 
    \begin{itemize}
    \item For each connected component $Q$ of $G[V_1]$ that is a GDP-tree, delete vertices of $Q$ that are not adjacent to any vertex of $V_2$ in $G-X$ and that are not cut-vertices of $Q-X$.
    \item For each connected component $Q$ of $G[V_1]$ that is a GDP-tree and for which $u_i$ is the last vertex in the ordering $u_1,u_2,\ldots,u_p$ having a neighbour in $Q$, choose a non-root vertex $w$ of a leaf block of $Q-X$. As $w$ is not a cut-vertex of $Q-X$, $w$ is adjacent to   $u_i$. Apply the operation $\ed_{(u_i,w)}$. 
    \item Delete vertex $u_i$.
\end{itemize}

It can be shown that all these operations are legal. Denote by $(G'',f'')$ the resulting pair after all  vertices of $V_2$ are processed (note that all vertices of $V_2$ are deleted).  

Assume that $Q$ is a connected component of $G[V_1]$ which is not
a GDP-tree. No vertex of $Q$ is deleted during the processing of
$V_2$. Hence, $f''(x)=d_Q(x)-1$ for every $x\in V(Q)$.

Since $Q$ is connected and is not a GDP-tree, Lemma~\ref{lem-degree}
implies that $Q$ is arc-weighted $(d_Q-1)$-degenerate. Therefore,
$Q$ is arc-weighted $f''|_{V(Q)}$-degenerate. By
Lemma~\ref{lem-alternate}, there is a sequence of legal arc-deletion
operations which transforms $(Q,f''|_{V(Q)})$ into an edgeless pair $(Q_0,g)$ with $g\ge 0$. All vertices of $Q_0$ are isolated and hence can subsequently be deleted legally by
vertex-deletion operations.

Assume $Q$ is a connected component of $G[V_1]$ and a GDP-tree. The graph $Q-X$ is connected and contains a vertex $w$ such that $f''(w)\ge d_{Q-X}(w)$, while $f''(x)\ge d_{Q-X}(x)-1$ for every other vertex $x\in V(Q-X)$. By repeatedly deleting a leaf of a spanning tree of $Q-X$ rooted at
$w$, and deleting $w$ last, all vertices can be deleted legally, as
the required inequalities are preserved after each deletion. 
\end{proof}

 \section{Some remarks and questions}

As a generalization of degeneracy, Bernshteyn and Lee \cite{BL} introduced the concept of weak degeneracy. 
 
A graph  $G$ is  $f$-degenerate if starting with the pair $(G,f)$,
we can apply legal vertex-deletion operations recursively to delete all vertices. Assume that each vertex $x$ has a set $L(x)$ of $f(x)+1$ acceptable colours. The vertex-deletion operation $\dd_x$ is interpreted as $x$ being coloured and each of its neighbours $y$ losing one permissible colour. Thus applying $\dd_x$  to $(G,f)$ results in $(G-x, f')$, where $f' = f-1_{N_G(x)}$. The   operation $\dd_x$ is legal if $f' \ge 0$.  

For the definition of weak degeneracy, another operation is used, called delete-save operation. This operation is defined as follows:

		\begin{definition} 
		Assume $(G,f) \in \mathcal{L}$, $x \in V(G)$ and $y \in N_G(x)$. The {\em delete-save operation} $\ds_{(x,y)}$ on $(G,f)$ is defined as   
			$${\ds_{(x,y)}(G,f)=(G-x,f_{-x}+1_y). }$$ 
			The   operation $\ds_{(x,y)}$ is    {\em legal} for $(G,f)$ if $f_{-x} +1_y  \ge 0$  and  $f(x) > f(y)$.  
		\end{definition}

Assume that $L$ is a list assignment of $G$ with $|L(v)| = f(v)+1$ for each vertex $v$. The delete-save operation $\ds_{(x,y)}$ means that we colour vertex $x$ with a colour in $L(x) - L(y)$, and hence $x$ is deleted, and each uncoloured neighbour of $x$ loses one permissible colour, however,  $y$ does not lose any colour. For $x$ to be coloured in this way, the set $L(x) - L(y)$ needs to be non-empty. Thus, the operation is legal if $f(x) > f(y)$.

We say $G$ is \emph{weak $f$-degenerate} if there is a sequence $\Omega$ of legal operations such that $\Omega(G,f) = (\emptyset, f')$, where $\emptyset$ is a non-vertex graph and $f' \ge 0$, where each operation is either a vertex-deletion operation or a delete-save operation.

Observe that $\ds_{(x,y)}$ applied to $(G,f)$ has the same effect as $\dd_x \ed_{(x,y)}$ applied to $(G,f)$. Therefore, if $G$ is weak $f$-degenerate, then it is arc-weighted  $f$-degenerate. The following lemma shows that weak  $f$-degeneracy corresponds to a special type of arc-weighted acyclic orientation. The proof of this lemma is straightforward and has been omitted.

\begin{lem}
    \label{lem-weak}
    A graph $G$ is weak  $f$-degenerate if and only if there is an arc-weighted orientation $(D,w)$ such that $d_{(D,w)}^+(v) \le f(v)$ for each vertex $v$, and for each  sub-digraph $D'$ of $D$,  either $D'$ has a sink, or $(D', w)$ has a dominating arc $(u,v)$ such that $u$ is a sink in $D'-(u,v)$. 
\end{lem}

Assume $L$ is an $(f+1)$-list assignment of $G$.  If $G$ is $f$-degenerate (respectively, weak $f$-degenerate or arc-weighted $f$-degenerate), then there is an algorithm, which we call the D-Algorithm (respectively, the WD-Algorithm or AWD-algorithm), that finds an $L$-colouring of $G$.

Both the D-Algorithm and the WD-Algorithm delete  one vertex in each move. The difference is that in the WD-Algorithm, a move that deletes vertex $v$ may "protect" one neighbour of $v$ from losing a colour. 
  
  The AWD-Algorithm  deletes one edge in each move. 
The advantage of the  operation $\ed_{(x,y),s}$ lies in its flexibility: it prohibits the vertex $x$ from using any colour assigned to the vertex $y$ without immediately colouring $x$. 
This flexibility is  useful in many cases, for example in Thomassen's proof of the 5-choosability of planar graphs. In \cite{BL}, a proof was given for the weak $4$-degeneracy of planar graphs. However, that proof has a gap. A correct and more complicated proof of this result was given in \cite{BLS2024}. Our paper has gone through a series of revisions. The first version of this paper appeared on arXiv in August 2023. At that time, there was no correct proof of the weak $4$-degeneracy of planar graphs. Also, the proof of arc-weighted $4$-degeneracy of planar graphs is   shorter and more transparent than the correct proof of their weak $4$-degeneracy.

The weak  degeneracy $d^*(G)$  of a graph $G$ is the minimum $k$ such that $G$ is weak   $k$-degenerate. From the definition it follows that $ d_{w}(G) \le d^*(G) $. In \cite{Yang} it was proved that for any integer $r$, there are $r$-regular graphs $G$ with $d^*(G) \le \lfloor \frac r2 \rfloor + 2$. As $r$-regular graphs have  degeneracy $r$, we know that the difference $d(G)- d^*(G)$ can be arbitrarily large, and hence the difference 
$d(G)-d_{w}(G) $ can be arbitrarily large.  The difference between $d^*(G)$ and $d_w(G)$ can also be arbitrarily large. The argument of Bradshaw and Zeng~\cite{BZ25} can be interpreted in our framework to give, for some absolute constant $c>0$, $d_w(K_{n,n}) \le n-c\sqrt{n\log n}$. On the other hand, Bernshteyn and Lee~\cite{BL} showed that
$d^*(K_{n,n}) \ge n-\sqrt{2n}$. Consequently, $$d^*(K_{n,n})-d_w(K_{n,n})
\ge c\sqrt{n\log n}-\sqrt{2n},$$
 which tends to infinity as $n\to\infty$.

In earlier versions of this paper, different generalizations of acyclic orientation are studied, and  the arc-weighted acyclic orientation $(D,w)$ in the current paper was called Type-2 acyclic in one of the earlier versions. 
The Type-1 acyclic arc-weighted orientation is defined as follows:

\begin{definition}
    \label{def-type2}
 We say that an arc-weighted digraph $(D,w)$ is \emph{Type-1} acyclic if $(D,w)$ contains no  nonempty Eulerian sub-digraph, i.e. a sub-digraph $(D', w)$ with $d_{(D',w)}^+(v) = d_{(D', w)}^-(v)$ for each vertex $v$.
\end{definition}
 
The Type-1 acyclic arc-weighted digraph $(D,w)$ is also a natural generalization of the acyclic digraph to the arc-weighted acyclic digraph. Compared to the Type-2 arc-weighted acyclic digraph (which we simply call arc-weighted acyclic digraph in this paper), there are two undesirable features of Type-1 acyclic arc-weighted digraphs:

(1) It is co-NP-complete to determine whether a given arc-weighted digraph is Type-1 acyclic. Equivalently, determining whether a given arc-weighted digraph contains a non-empty Eulerian sub-digraph is NP-complete.  This is through a reduction from the zero-sum subset problem: given a set $S$ of integers, the zero-sum subset problem asks if $S$ has a non-empty subset $S'$ with $\Sigma_{a \in S'}a=0$. It is well-known that zero-sum subset problem is NP-complete.
Let $S^+=\{a \in S: a > 0\}$ and $S^-=\{a \in S: a < 0\}$.  Let $ (D, w)$ be the arc-weighted digraph consisting of two vertices $u,v$,  $|S^+|$ parallel arcs $(u,v)$ of weights $S^+$ and $|S^-|$ parallel arcs $(v, u)$ of weights $\{|a|: a \in S^-\}$. Then $D$ has a non-empty sub-digraph $D'$ with $d_{(D',w)}^+(u) = d_{(D',w)}^-(u)$  if and only if $S$ has a non-empty subset $S'$ with sum $\Sigma_{a \in S'}a=0$. 

(2) Given a Type-1 acyclic arc-weighted digraph $(D,w)$, let $f(v)= d_{(D, w)}^+(v)+1$ for each vertex $v$, we know that $G$ is $f$-AT and therefore $f$-paintable and $f$-choosable. However, given an $f$-list assignment $L$ of $G$, we do not know of an efficient algorithm to find an $L$-colouring of $G$.   It is  an interesting and  challenging problem to find an efficient algorithm (even a sub-exponential algorithm) that constructs an $L$-colouring from such an orientation.

Also, the following question remains open.

\begin{question}
    \label{q2}
    Given a Type-1 acyclic arc-weighted digraph $(D,w)$, let $f(v)= d_{(D, w)}^+(v)+1$ for each vertex $v$. Is $G$  $f$-DP-colourable? Or even $f$-DP-paintable?
\end{question} 

As shown in Lemma \ref{lem-weak}, weak degeneracy corresponds to a special family of arc-weighted acyclic orientation $(D,w)$. This was called Type-4 acyclic orientation in an earlier version of this paper. Type-3 acyclic orientation corresponds to a generalization of weak degeneracy, where the operation 
$\ds_{(x,y)}$ is replaced by $\ds_{(x,U)}$ where $U$ is a subset of $N_G(x)$: When colouring a vertex $x$, instead of protecting one neighbour $y$ of $x$ (by colouring $x$ with a colour from $L(x)-L(y)$), we may protect a subset $U$ of $N_G(x)$ (by colouring $x$ with a colour from $L(x)- \cup_{y \in U}L(y)$, provided that $|L(x)| > \sum_{y \in U}|L(y)|$). In this language, the vertex deletion operation $\dd_x$ is the same as $\ds_{x, \emptyset}$. For the reasons mentioned above and for brevity, we have decided to concentrate on Type-2 arc-weighted acyclic orientations in the current version of the paper. However, the other types of arc-weighted acyclic orientation are also useful and have been studied in some published papers \cite{BZ25,XZ25, LWY}, where the earlier version of this paper is cited. 

For a given graph $G$ and $f \in \mathbb{N}^G$, it can be determined in  linear time whether $G$ is $f$-degenerate. Given an arc-weighted orientation  $(D,w)$ of a graph $G$, a naive peeling procedure determines in  time $O(m^2)$  whether $(D,w)$ is acyclic. The  following question remains open:

\begin{question}
    \label{q3}
     What is the complexity of the following problem:
     
     Instance: A graph $G$ and a function $f \in \mathbb{N}^G$.
     
     Question: Is $G$ arc-weighted  $f$-degenerate?  
\end{question}

  \section*{Acknowledgment}
  The authors thank Patrice Ossona De Mendez, Jakub Kozik, Grzegorz Gutowski, and Reza Naserasr for valuable discussions. We thank anonymous reviewers for their careful reading of the manuscript and for many valuable comments that improved the presentation of this paper.

\bibliographystyle{plain}
\bibliography{template}

	\end{document}